\documentclass[a4paper]{amsart}\usepackage{amsfonts,amssymb,bbm,graphicx,enumerate,rotating}

\def\C{\Bbb{C}}\def\k{\mathbbm{k}}
\def\N{\Bbb{N}}
\def\R{\Bbb{R}}\def\Z{\Bbb{Z}}

\def\li{~\\ $\bullet$ }\newcommand{\ls}{~\\ $\star$ }\def\di{\partial}

\def\suml{\sum\limits}
\def\liml{\lim\limits}

\def\capl{\mathop\cap\limits}

\newcommand{\quotient}[2]{{\left.\raisebox{1.6ex}{$#1$}\!\!\!\!\!{\scalebox{2}{\ensuremath\diagup}}
\!\!\!\!\!\raisebox{-1ex}{$#2$}\right.}}
\newcommand{\quotients}[2]{{\footnotesize\left.\raisebox{0.4ex}{$#1$}\! / \!\raisebox{-0.4ex}{$#2$}\right.}}

\def\tu{\tilde{u}}\def\tv{\tilde{v}}

\def\ha{{\hat{a}}}\def\he{\hat{e}}
\def\hF{\hat{F}}\def\hG{{\widehat{G}}}\def\hg{{\hat{g}}}
\def\hLa{\widehat{\La}}
\def\hM{{\widehat{M}}}\def\hR{{\widehat{R}}}\def\hU{\hat{U}}
\def\hz{\hat{z}}

\def\al{\alpha}\def\ga{\gamma}\def\be{\beta}\def\De{\Delta}
\def\ep{\epsilon}
\def\la{\lambda}\def\La{\Lambda}\def\si{\sigma}\def\Si{\Sigma}

\def\cD{\mathcal D}
\def\cG{\mathcal G}\def\cK{{\mathcal K}}
\def\cm{{\frak m\hspace{0.05cm}}}\def\cR{\mathcal{R}}

\def\ua{{\underline{a}}}\def\ub{{\underline{b}}}
\def\uF{\underline{F}}\def\uG{\underline{G}}\def\uH{\underline{H}}

\def\uv{\underline{v}}\def\uw{{\underline{w}}}\def\ux{\underline{x}}\def\uy{{\underline{y}}}\def\uz{{\underline{z}}}

\def\one{{1\hspace{-0.1cm}\rm I}}\def\zero{\mathbb{O}}

\newcommand{\ber}{\begin{array}{l}}\newcommand{\eer}{\end{array}}
\newcommand{\bpm}{\begin{pmatrix}}\newcommand{\epm}{\end{pmatrix}}
\newcommand{\bM}{\begin{matrix}}\newcommand{\eM}{\end{matrix}}
\newcommand{\bee}{\begin{enumerate}}\newcommand{\eee}{\end{enumerate}}
\newcommand{\bei}{\begin{itemize}}\newcommand{\eei}{\end{itemize}}

\def\wrt{with respect to }\def\sset{\subset}\def\sseteq{\subseteq}\def\ssetneq{\subsetneq}\def\smin{\setminus}
\def\jbj{jet-by-jet}
\def\iff{if and only if\ }
\def\Mat{Mat(m,n;R)}

\def\bull{\vrule height .9ex width .9ex depth -.1ex }

\newtheorem{Lemma}{Lemma}[section]\newcommand{\bel}{\begin{Lemma}}\newcommand{\eel}{\end{Lemma}}
\newtheorem{Theorem}[Lemma]{Theorem}\newcommand{\bthe}{\begin{Theorem}}\newcommand{\ethe}{\end{Theorem}}
\newtheorem{Proposition}[Lemma]{Proposition}\newcommand{\bprop}{\begin{Proposition}}\newcommand{\eprop}{\end{Proposition}}
\newtheorem{Corollary}[Lemma]{Corollary}\newcommand{\bcor}{\begin{Corollary}}\newcommand{\ecor}{\end{Corollary}}
\newtheorem{Definition}[Lemma]{Definition}\newcommand{\bed}{\begin{Definition}}\newcommand{\eed}{\end{Definition}}
\newtheorem{Definition-Proposition}[Lemma]{Definition-Proposition}

\def\bpr{~\\{\em Proof.\ }}\newcommand{\epr}{$\bull$\\}
\newtheorem{Remark}[Lemma]{Remark}\newcommand{\beR}{\begin{Remark}\rm}\newcommand{\eeR}{\end{Remark}}
\newtheorem{Example}[Lemma]{Example}\newcommand{\bex}{\begin{Example}\rm}\newcommand{\eex}{\end{Example}}
\newtheorem{Problem}[Lemma]{Problem}\newcommand{\bprob}{\begin{Problem}\rm}\newcommand{\eprob}{\end{Problem}}

\newcommand{\bet}{\begin{tabular}{cccccccc}}\newcommand{\eet}{\end{tabular}}
\newcommand{\beq}{\begin{equation}}\newcommand{\eeq}{\end{equation}}

\def\kpd{{$\k$-polynomially-defined}}

\newcommand\isom{\xrightarrow{\,\smash{\raisebox{-0.65ex}{\ensuremath{\sim}}}\,}}

\title[]{G\MakeLowercase{roup actions on filtered modules and finite determinacy.}\\F\MakeLowercase{inding
large submodules in the orbit by linearization}}
\author[]{G\MakeLowercase{enrich} B\MakeLowercase{elitskii and} D\MakeLowercase{mitry} K\MakeLowercase{erner}}
\address{Department of Mathematics, Ben Gurion University of the Negev, P.O.B. 653, Be'er Sheva 84105, Israel.}
\email{genrich@math.bgu.ac.il}
\email{dmitry.kerner@gmail.com}
\date{\today}
\thanks{D.K. was partially supported by the postdoctoral fellowship
 at the University of Toronto and by the grant FP7-People-MCA-CIG, 334347.}

\begin{document}\maketitle
\begin{abstract}
Let $M$ be a module over a local ring $R$, with a group action $G\circlearrowright M$, not necessarily $R$-linear.
  To understand how large is the $G$-orbit of an element $z\in M$ one looks for the large submodules
of $M$ lying in $Gz$.
We provide the corresponding (necessary/sufficient) conditions in terms of the tangent space to the orbit,
  $T_{(Gz,z)}$.

 \

This question originates from the classical finite determinacy problem of Singularity Theory. Our treatment is
rather general, in particular we extend the classical criteria of Mather (and many others) to
 a broad class of rings, modules and group actions.

When a particular `deformation space' is prescribed, $\Si\sseteq M$,
  the determinacy question is translated into the properties of the tangent spaces, $T_{(Gz,z)}$, $T_{(\Si,z)}$, and in particular to the annihilator of their quotient, $ann\quotients{T_{(\Si,z)}}{T_{(Gz,z)}}$.
\end{abstract}\setcounter{secnumdepth}{6} \setcounter{tocdepth}{1}

\section{Introduction}

\subsection{Setup}
Let $R$ be a (commutative, associative) local ring over a base field $\k$ of zero characteristic.
Denote by $\cm\sset R$  the maximal ideal.
 (In the simplest case $R$ can be a {\em regular} ring, e.g. the rational functions regular at the origin, $\k[x_1,\dots,x_p]_{(\cm)}$; the formal power series, $\k[\![x_1\,\dots,x_p]\!]$; the converging power series, $\C\{x_1,\dots,x_p\}$; the smooth function germs, $C^\infty(\R^p,0)$.)
Geometrically, $R$ is the ring of regular functions on the (algebraic/formal/analytic etc.) germ $Spec(R)$.
We use some Artin-type approximation properties of $R$,  this excludes from consideration the rings like $C^r(\R^p,0)$, for $r<\infty$.

\

Let $M$ be a  module over $R$, with a descending filtration, $M=M_0\supsetneq M_1\supsetneq M_2\supsetneq\cdots$. This filtration defines the linear topology on $M$, the open neighborhoods of $z\in M$ are the sets $\big\{\{z\}+M_i\big\}_i$.
The simplest filtration is defined by the powers of an ideal $J$, namely $M_i=J^i\cdot M$.
More generally, any filtration by ideals, $R=J_0\supsetneq J_1\supsetneq J_2\supsetneq\cdots$, induces the corresponding filtration $M_i=J_i\cdot M$.

\

Fix a $\k$-linear group action $G\circlearrowright M$. (If a filtration $M_\bullet$ is given we usually assume the action filtered, i.e., $G\cdot M_i=M_i$.)  We address the classical question.
\beq
\text{\em For a `small' deformation $z\rightsquigarrow z'$, with $z,z'\in M$, are the initial and deformed elements  $G$-equivalent?}
\eeq
More precisely,  `does the orbit $Gz$ contain some open neighborhood $\{z\}+M_i$?'

\

In various applications one deforms $z$ not inside the whole module $M$, but inside some subspace $\Si\sseteq M$ of `allowed' deformations. We always assume that $\Si$ is $G$-invariant. Usually $\Si$ is `reasonably good', e.g.  is defined by some power series equations (or just by  polynomials).
 The topology on $\Si$ is induced from that on $M$, the open neighborhoods of $z\in\Si$ being of the form $\big(\{z\}+M_j\big)\cap \Si$.
\bex\label{Ex.Intro.Typical.Groups} For the cases below we fix
$\Si=M$.
\bee
\item Denote by $GL_R(M)$ the group of all the $R$-linear
automorphisms of $M$. For a module with filtration, $M_\bullet$,
denote by $GL^{(0)}_R(M)\sseteq GL_R(M)$ the subgroup of
automorphisms that preserve the filtration.
\item  Denote by $Aut_\k(R)$ the group of $\k$-linear automorphisms of the ring. Suppose $M$ is free and fix a set of generators, $\{e_i\}$, of $M$. Then $Aut_\k(R)$ acts on $M$, by $\sum_j a_j e_j\to \sum_j \phi(a_j)e_j$, for $a_j\in R$. This action depends essentially on the choice of $\{e_j\}$, but is well defined otherwise.
\item Suppose $M$ is free, of rank $mn$. Identify it with the space of $m\times n$ matrices over $R$, i.e., $M\isom \Mat$. Various subgroups of $GL_R(M)$ are related to the rich matrix structure. For example, the left multiplications $G_l:=GL(m,R)$,  the right multiplications $G_r:=GL(n,R)$,  the two-sided multiplications $G_{lr}:=G_l\times G_r$,
$A\to UAV^{-1}$.
\eee\eex
\bex\label{Ex.Intro.Typical.Def.Spaces}
Consider the module of square matrices, $M\isom Mat(m,m;R)$.
\bee
\item The congruences, $G_{congr}\approx GL(m,R)$, act by $A\to UAU^t$ and preserve the submodules of symmetric/anti-symmetric matrices. Therefore, for $A$ symmetric, it is natural to choose $\Si=Mat^{sym}(m,m;R)$, while  in the anti-symmetric case one chooses $\Si=Mat^{anti-sym}(m,m;R)$.
\item The conjugations, $G_{conj}\approx GL(m,R)$, act by $A\to UAU^{-1}$, and preserve the characteristic polynomial $det(\la\one-A)$. In this case one often chooses $\Si=\Si_A=\{B|\ det(\la\one-A)=det(\la\one-B)\}\sset Mat(m,m;R)$.
\eee\eex

\

\subsection{The finite and infinite determinacies.}\label{Sec.Intro.finite.determinacy}
Fix a module $M$ with filtration  $M_\bullet$, a group action $G\circlearrowright M_\bullet$ and a deformation subspace $\Si\sseteq M$. Suppose $\Si$ is $G$-invariant. An element $z\in\Si$ is called
$k$-$(\Si,G,M_\bullet)$-determined if $z\stackrel{G}{\sim}z'$ whenever $z'-z\in (\Si-\{z\})\cap M_{k+1}$.
More precisely:
\bed
The order of $(\Si,G,M_\bullet)$-determinacy is:
$ord_G^\Si(z)=min\Big\{k|\ Gz\supseteq (\{z\}+M_{k+1})\cap\Si\Big\}\le\infty$.
\eed
An element $z$ is  finitely-$G$-determined, that is $ord^\Si_G(z)<\infty$, if the orbit $Gz\sset\Si$ contains an open neighborhood of $z$ in the filtration topology. Finite determinacy means that $z$ is determined (up to $G$-equivalence) by its image in  $\quotients{M}{M_{k+1}}$ for some finite $k$.

Sometimes the filtration (eventhough strictly decreasing) contains `flat elements', i.e., $M_\infty=\capl_{i=0}^\infty M_i\neq\{0\}$. (The typical example is: $\{M_j=\cm^j\cdot M\}_j$ for the ring $R=C^\infty(\R^p,0)$, so $\cm^\infty\neq\{0\}$.)
 An element $z$ is called infinitely-$(\Si,G)$-determined  if
 $z'\stackrel{G}{\sim}z$ whenever $z'-z\in (\Si-\{z\})\cap M_\infty$. The infinite determinacy means that $z$
 is determined by its image in the completion $\hM$ of $M$ \wrt $M_\bullet$.
\bex
For the filtration $M_i=J^i\cdot M$ take the completion of $M$ \wrt $M_\bullet$. Any element $z\in M$ maps to $\hz\in\hM$ that is presentable as a power series in the generators of $J$. This power series is the ``Taylor expansion at the origin". Finite determinacy means that $z$ is fixed (up to the $G$-action) by a finite number of terms in its Taylor expansion. Infinite determinacy means that the  ``full Taylor expansion" fixes $z$ up to the $G$-action.
\eex
\bex
Suppose $M$ is a free $R$-module, identify it with $R^{\oplus n}$. 
 Suppose $R$ is one of the classical local `geometric' rings, $\k[[\ux]]$ or $\k\{\ux\}$, when $\k$ is a normed field.
\bee\item If $n=1$ then $z\in \cm\cdot M=\cm$ defines a (formal/analytic) hypersurface singularity at the origin, $\{z=0\}\sset (\k^p,0)$. The group $Aut_\k(R)$ then coincides with the group of the local coordinate-changes, $\cR$. Thus we get the classical {\em right-equivalence} of the Singularity Theory. Similarly, the group $GL(1,R)\rtimes Aut_\k(R)$ induces the classical {\em contact} equivalence, $\cK$.
\item Suppose $n>1$, so that $z$ is an $n$-tuple in $R^{\oplus n}$. Assume all the entries of $z$ belong to $\cm$, i.e., ``$z$ vanishes at the origin of $Spec(R)$". Then $z$ can be considered as a map from $Spec(R)$ to $(\k^n,0)$.
 Again we get the classically studied equivalences of maps, the right, $Aut_\k(R)=\cR$, and the contact, $GL(n,R)\rtimes Aut_\k(R)=\cK$.
\eee\eex
The finite determinacy of maps (or of the corresponding singularities) under various equivalences has been intensively studied since the seminal works  \cite{Mather1968}, \cite{Arnol'd}, \cite{Tougeron1968}.

\

For various group actions, $G\circlearrowright M$, the quotient  $\quotients{M}{G}$ parameterizes the geometric/algebraic objects and the determinacy bears important information about their deformation theory.
\bex
Continuing examples \ref{Ex.Intro.Typical.Groups} and \ref{Ex.Intro.Typical.Def.Spaces}.
\bee
\item Thinking of $A\in\Mat$ as a presentation matrix of the module $coker(A)$, in the projective resolution, we get:
   matrices up to the $G_{lr}$-equivalence correspond to the modules over $R$. Similarly, thinking of $A$ as the matrix of generators of $Im(A)$ we get:    matrices up to the $G_{r}$-equivalence correspond to the submodules of  $R^{\oplus m}$.
\item  The (anti-)symmetric matrix, $A\in Mat(m,m;R)$, $A^t=\pm A$, considered up to the congruence, $A\stackrel{G_{congr}}{\sim}UAU^t$, defines a (skew-)symmetric form over $R$.
\item A quadratic matrix considered up to the conjugation, $A\stackrel{G_{conj}}{\sim}UAU^{-1}$ corresponds to a representation (of a group/algebra/etc.).
\eee\eex
Finite determinacy in these cases implies that the deformation theory is essentially finite dimensional.

\subsection{}
The question of determinacy can be restated as:
\beq
\text{which deformations of $z$ inside $\Si$ are irrelevant, i.e., lie inside the orbit $Gz$?}
\eeq
In other words: {\em how large is the orbit $Gz$ as compared to $\Si$}? In this paper we linearize this question, i.e., transform it to comparison of the tangent spaces, $T_{(Gz,z)}\sseteq T_{(\Si,z)}$.  We prove
the Mather-type determinacy criterion in great generality, for a large class of rings, group-actions and deformation spaces.
This reduces the determinacy problem of $(\Si,G,z)$ to the study of the quotient of tangent spaces, $\quotients{T_{(\Si,z)}}{T_{(Gz,z)}}$.

\subsection{Contents and the structure of the paper}
The main results are formulated in section \ref{Sec.Results}. In section \ref{Sec.Intro.Relation.to.Singular.Theory} we give a brief historical sketch and relate our work to the vast field of results on finite determinacy in Singularity Theory.

\

We work with a broad class of groups and to our knowledge in this generality the tangent space $T_{(G,g)}$ to $G$ at $g\in G$ has not yet been defined. Therefore in section \ref{Sec.Group.Actions} we lay some foundations.

Fix an action $G\circlearrowright M$, a filtration on $M$ induces the filtration on $G$.
(For example, $G^{(0)}\sseteq G$ is the subgroup of those elements that preserve the filtration; $G^{(1)}\sseteq G^{(0)}$ is the subgroup of those elements that act unipotently \wrt the filtration.)
Lemma \ref{Thm.Group.Actions.Completion.General} summarizes the functorial properties of this filtration, the related projections and completions. Then we specify the class of groups we work with. As is seen from examples of the introduction, this class must contain various subgroups of $GL_R(M)\rtimes Aut_\k(R)$. As the local ring $R$ is not necessarily $\quotients{\k[[\ux]]}{I}$ (or a subring of this), it is not natural to restrict to the subgroups defined by some power series equations over $R$. Rather, we consider all the $\k$-linear (not necessarily $R$-linear!) endomorphisms, $End_\k(M)$, and the $\k$-linear automorphisms $GL_\k(M)$. Note that $M$ is uncountably generated as a $\k$-vector space, therefore $End_\k(M)$, $GL_\k(M)$ are huge. We call a subgroup $G\sseteq GL_\k(M)$ {\em ``\kpd"} if its defining conditions translate into a system of polynomial equations in some (and hence any) Hamel basis.
 (This system is usually uncountable and involves uncountable number of variables.) The class of \kpd-groups is very large (lemma \ref{Thm.Group.Actions.kpd.groups}), in particular it is much larger than the class of groups algebraic over $R$ or projective limits of polynomially defined (over $R$) or those defined by formal power series.

As the germ $(G,g)$ is algebraic {\em over $\k$} (though of uncountable dimension and codimension)
the tangent space  is defined as the Zariski tangent space, $T_{(G,g)}$. This definition of $T_{(G,g)}$ goes via the embedding $G\sseteq GL_\k(M)$,
i.e., is not internal, and various pathologies might occur. To prevent this we restrict to a subclass of \kpd-groups: the {\em groups of Lie type},
  section \ref{Sec.Group.Actions.Formally.Smooth.Groups}.
  These are the groups admitting some substitution of the logarithmic/exponential maps.

Finally, when both $M$ and $G$ are complete  and moreover $G$ is unipotent \wrt the filtration we give an alternative definition of $T_{(G,\one)}$,
in the ``internal/canonical" way: via the logarithm/exponential maps, section \ref{Sec.Group.Actions.Logarithm.Exponent}. If $G$ is \kpd\ and of Lie type,
 then the two definitions coincide and $G$ becomes a Lie group (proposition \ref{Thm.Group.Action.Tangent.Space.Comparison.Lie.Group}).
Thus part of section \ref{Sec.Group.Actions} can be thought as the axiomatization of the classical Lie-group notion.
In section \ref{Sec.Groups.Acting.on.Matrices} we consider the simplest examples, $M=\Mat$ and $G\sseteq G_{lr}\rtimes Aut_\k(R)$.

\

Our main result is theorem  \ref{Thm.Intro.Linearization}. The proof goes in three steps.
\bei
\item In section \ref{Sec.Linearization.procedure.jbj} we define the determinacy on locally filtered (not necessarily linear) sets. We
pass from the tangent space condition, $z'-z\in T_{(Gz,z)}$,
  to the \jbj-equivalence, $z'\stackrel{G_{jbj}}{\sim}z$. The later means: for any $q$ holds $z'\in Gz+M_q$. This transition is done in the general
  setting: $M$ is a filtered $\k$-vector space, the action $G\circlearrowright M$ is filtered, and $M,G$ are complete
  \wrt\ filtrations. (Here $G$ is not necessarily \kpd\ and $M$ is of arbitrary dimension.)
\item One turns the \jbj-equivalence into the equivalence of $\cm$-adic completions, $\hz'\stackrel{\hG}{\sim}\hz\in\hM^{(\cm)}$.
 (Here we use Popescu's theorem, stated in section \ref{Sec.Approximations.Popescu.Thm}.)
\item Eventually one passes  from the equivalence of completions, $\hz'\stackrel{\hG}{\sim}\hz$, to the ordinary equivalence, $z'\stackrel{G}{\sim}z$. (Here we use the Artin-type approximation theorems, stated in section \ref{Sec.Approximations.Artin}.)
\eei
While the first step is done without many assumptions, for the
second and third steps we impose some restrictions: the module $M$
is finitely generated over $R$ and the group action
$G\circlearrowright M$ is good enough (\kpd\ and $G^{(1)}\sset G$ is of Lie type). We use these
assumptions in the proofs, but we hope they can be significantly
weakened (in the future) by using some stronger approximation
results.

All these results are combined in section \ref{Sec.Linearization.Final} to finish the proof of theorem \ref{Thm.Intro.Linearization}. Then we prove the corresponding criterion for finite determinacy in families.

\subsection{Acknowledgements}

Many thanks are to V.Goryunov, G.M.Greuel, T.H.Pham and I.Tyomkin.
The extraordinary work of the referee not only helped to improve the exposition but forced us to generalize many statements and to simplify their proofs.

\section{The main results and remarks}\label{Sec.Results}

\subsection{Notations}
We use the space $End_\k(M)$ of all the $\k$-linear endomorphisms, here $M$ is considered just as a $\k$-vector space. Accordingly $GL_\k(M)$ is the group of all the $\k$-linear automorphisms.

 We denote the zero matrix (or the zero element of $End_\k(M)$) by $\zero$, while the identity matrix/endomorphism by $\one$.
By $\one$ or $\one_G$ we denote also the unit element of the group $G$.

The group $G$ is always a subgroup of $GL_\k(M)$, i.e., $G$ comes with its (faithful) action $G\circlearrowright M$.

The {\em unipotent} subgroup $G^{(1)}\sseteq G$ is defined by
\beq\label{Eq.def.of.G1}
G^{(1)}:=\{g\in G|\ \forall j:\ [g]=[Id]\circlearrowright \quotient{M_j}{M_{j+1}}\}.
\eeq
For example, suppose $M=R^{\oplus n}$ and the filtration is $M_i=\cm^i\cdot M$. Then $GL_R(M)$ is the group of all the $n\times n$ matrices invertible over $R$, while $GL^{(1)}_R(M)=\{\one+U|\ U\in Mat(n,n;\cm)\}$.

\subsection{Assumptions}
Though $R$ is not necessarily Noetherian, we assume that the $\cm$-adic completion, $\hR$, is Noetherian. Thus  while we allow rings like $C^\infty(\R^p,0)$, we do not allow rings like $\k[[x_1,x_2,\dots]]$.

We always assume that the $R$-module $M$ is finitely generated (over $R$).
In the filtered case we assume that all $M_j$ are finitely generated.
We always assume that the filtration is `essentially decreasing', i.e., satisfies $\capl_j M_j\sseteq\capl_i(\cm^i\cdot M)$. Equivalently: for any $i$ there exists $k_i$ such that $M_{k_i}\sseteq\cm^i\cdot M$.

In most sections we assume that
 the action $G\circlearrowright M$ is good enough, namely:
\beq\label{Eq.assumptions.of.kpd.fs.etc}\ber
\text{\em the subgroup $G\sset GL_\k(M)$ and its completion $\widehat{G}\sseteq GL_\k(\hM)$ are \kpd;}\\
\text{\em their unipotent parts, $G^{(1)}$, $\widehat{G}^{(1)}$, are of Lie type}.
\eer\eeq
These conditions are stated/studied in section \ref{Sec.Group.Actions.kpd.groups}, section \ref{Sec.Group.Actions.Tangent.Space}, section \ref{Sec.Group.Actions.Formally.Smooth.Groups}.
In this generality we define the tangent spaces, $T_{(Gz,z)}\sseteq T_{(\Si,z)}\sseteq T_{(M,z)}$. Initially these are just $\k$-vector subspaces but in some places we assume that they are $R$-submodules of $T_{(M,z)}$, this imposes some restrictions on $G$ and $\Si$.

The condition  ``$(G,\one)$ is of Lie type" ensures that the tangent space $T_{(Gz,z)}$ ``approximates" the germ $(Gz,z)$.

\

In many statements we use the assumption ``$R$ has the relevant approximation property". This condition depends on the type of equations, see section \ref{Sec.Approximation.Propert.Groups}, e.g. the Artin approximation property of $R$ suffices for polynomial (or analytic for $R=\C\{\ux\}$) equations.

\subsection{Transition to the tangent spaces}
As one sees in example \ref{Ex.Intro.Typical.Def.Spaces} the deformation space $\Si\sseteq M$ can be a highly non-linear subset. Still, the finite determinacy implies that some projections of $\Si$ contain large linear subspaces.
\bel\label{Thm.Intro.Linear.is.necessary} Fix a filtration $M_\bullet$ and a filtered action $G^{(1)}\circlearrowright M_\bullet$.
Suppose $z$ is $r$-($\Si,G^{(1)}$)-determined. Then for any $i\ge r$: $\Big(\Si-\{z\}\Big)\cap M_i+M_{i+1}$ is an $R$-submodule of $M$.
\vspace{-0.2cm}\eel
 The proof is in section \ref{Sec.Linearization.jbj.Determinacy.Action.Filtered.Vector.Space}.

\

Therefore the finite determinacy depends on a more general question:
\beq\label{Eq.Linearized.Question}
\text{Find the largest submodule, $\La\sset M$, satisfying: $\{z\}+\La\sseteq Gz$.}
\eeq

\

As $M$ is an $R$-module, we identify  $T_{(M,z)}\approx M$, as $R$-modules.
 Accordingly we identify  $T_{(Gz,z)}$ with its image in $M$.
Our main result reduces the `linearized' question of equation \eqref{Eq.Linearized.Question} to the tangent space.

\bthe\label{Thm.Intro.Linearization}
Suppose the (filtered) action $G\circlearrowright\{M_i\}$ satisfies assumptions \eqref{Eq.assumptions.of.kpd.fs.etc}.  Suppose that $G$ is unipotent for the filtration $\{M_i\}$.
\\1. If $M_i\sseteq T_{(Gz,z)}$ and $R$ has the relevant approximation property then $\{z\}+M_i\sseteq Gz$.
\\2. Suppose $T_{(Gz,z)}\sseteq T_{(M,z)}$ is an $R$-submodule. If $\{z\}+M_i\sseteq Gz$ then $M_i\sseteq T_{(Gz,z)}$.
\ethe
\noindent(In part (2.) the assumption ``$T_{(Gz,z)}\sseteq T_{(M,z)}$ is an $R$-submodule" can be weakened to: ``$T_{(Gz,z)}\sseteq T_{(M,z)}\approx M$ is eventually a submodule" in the following sense:
 for some $N<\infty$ the intersection $T_{(Gz,z)}\cap M_N$ is a submodule of $M$, see remark \ref{Ex.Tangent.Space.Dont.need.submodule}.)

\beR\label{Ex.Counterexample.to.strengthening.1}
One might wish to strengthen part 1 to the statement: ``If $G$ is unipotent for some filtration on $M$ then $\{z\}+T_{(Gz,z)}\sseteq Gz$".
But in general $T_{(Gz,z)}$ is not invariant under the $G$ action and this inclusion does not hold. For example, let $M_i=Mat(m,n;\cm^i)$ and consider the action $G^{(1)}_{lr}\circlearrowright M$ of example \ref{Ex.Intro.Typical.Groups}. Here $T_{(G^{(1)}A,A)}$ is spanned by the matrices $\tu A$ and $A\tv$ for all the possible $(\tu,\tv)\in Mat(m,m;\cm)\times Mat(n,n;\cm)$, see section \ref{Sec.Groups.Acting.on.Matrices}. The inclusion $\{z\}+T_{(Gz,z)}\sseteq Gz$ would mean the solvability of the  equation $(\one+u)A(\one+v)=A+\tu A+A\tv$ for the  arbitrary choice
$(\tu,\tv)\in Mat(m,m;\cm)\times Mat(n,n;\cm)$  and the unknowns $(u,v)\in Mat(m,m;\cm)\times  Mat(n,n;\cm)$.
Take e.g. $A=\bpm 0&1\\0&0\epm$ and $\tu=\bpm 0&0\\x&0\epm=\tv$. By the direct check, in this case the system has no solutions. Note that $T_{(G^{(1)}A,A)}=Span_\cm\Big(\bpm 1 &0\\0&0\epm,\bpm 0 &1\\0&0\epm,\bpm 0 &0\\1&0\epm \Big)$ is not $G_{lr}$-invariant.
\eeR

\subsubsection{}
The theorem addresses the `pointwise' determinacy: one checks the equivalence $z'\stackrel{G}{\sim}z$ for a particular $z'$. The natural question is: when a family $\{z_t\}_{t\in(\k,0)}$ can be trivialized?
 (Namely, when there exists a family $\{g_t\}_{t\in(\k,0)}\sset G$ such that $\{z_t=g_t z\}_{t\in(\k,0)}$?)
We prove the corresponding criterion in section \ref{Sec.Determinacy.in.Families}.

\subsubsection{}
Usually we start with just an action $G\circlearrowright M$, with no prescribed filtration. To achieve the best criteria/bounds one looks for a suitable (optimal) filtration, thus one compares $T_{(Gz,z)}$ to $M$.  Two cases are possible:
\bee[i.]
\item
$T_{(Gz,z)}\supseteq J\cdot M$, for some ideal $\{0\}\ssetneq J\sset R$. In this case it is useful to consider the filtration $M_i=(\sqrt{J})^i\cdot M$. Here we take the radical, $\sqrt{J}$, to get the ``"most refined"' filtration and accordingly the biggest $G^{(1)}$, see section \ref{Sec.Intro.Annihilator}.
\item
$T_{(Gz,z)}\not\supset J\cdot M$, for any ideal $\{0\}\ssetneq J\sset R$, equivalently $rank_R(T_{(Gz,z)})< rank_R(M)$. Then one looks for the (biggest/simplest) submodule $\La\sset M$ that satisfies: $\La$ is $G$-invariant, $z\in\La$, and $T_{(Gz,z)}\supseteq J\cdot \La$, for some $\{0\}\ssetneq J\sset R$. For such a submodule one takes the restriction, $G|_\La$, i.e., considers $\La$ as the ambient module (rather than $M$). For most cases the subgroup $G|_\La\sseteq GL_\k(\La)$ again satisfies the assumption \eqref{Eq.assumptions.of.kpd.fs.etc}, hence one can use the theorem.
\eee

\subsection{The annihilator of the quotient and the order of determinacy}\label{Sec.Intro.Annihilator}
Usually  one begins with the action $G\circlearrowright \Si\sseteq M$, but with no prescribed filtration and no prescribed subgroup $G^{(1)}\sseteq G$. Theorem \ref{Thm.Intro.Linearization} gives a separate statement for each filtration of $M$ (or of $\La$). Among all these versions one would like to choose an optimal bound.
Using the embedding $T_{(Gz,z)}\sseteq T_{(\Si,z)}\sseteq T_{(M,z)}$ we get: the largest submodule of equation \eqref{Eq.Linearized.Question} satisfies: $\La\sseteq T_{(\Si,z)}$.
This translates the initial question into ``how large is $T_{(\Si,z)}$ as compared to $T_{(Gz,z)}$?".
To quantify this one usually studies their conductor, i.e., the annihilator of the quotient of the two modules:
\beq
ann\quotient{T_{(\Si,z)}}{T_{(Gz,z)}}:=\{f|\ f\cdot T_{(\Si,z)}\sseteq T_{(Gz,z)}\}\sset R
\eeq
In many cases of interest both tangent spaces are $R$-modules, and their quotient is an $R$-module as well. Then the annihilator is an ideal in $R$.

This annihilator shows how far is $T_{(Gz,z)}$ from $T_{(\Si,z)}$. It is defined via the tangent spaces and thus  controls the ``infinitesimal determinacy". By Theorem \ref{Thm.Intro.Linearization} this annihilator is tightly related to the standard determinacy:
\bcor\label{Thm.Intro.Reduction.to.annihilator}
Suppose $T_{(\Si,z)}\sset T_{(M,z)}$ is a finitely generated submodule and for a (finitely-generated) ideal $J\ssetneq R$
the filtration $\{J^i\cdot T_{(\Si,z)}\}_i$ is $G$-invariant. Suppose  the corresponding unipotent subgroup,  $G^{(1)}\sseteq G$,
 satisfies the assumption \eqref{Eq.assumptions.of.kpd.fs.etc}.
\\1. Suppose $R$ has the relevant approximation property and $J\sseteq ann\quotients{T_{(\Si,z)}}{T_{(Gz,z)}}$. Then
$\{z\}+ J\cdot \sqrt{J}\cdot T_{(\Si,z)}\sseteq Gz$.
\\1'.  If in addition $J\cdot T_{(\Si,z)}\sseteq T_{(G^{(1)}z,z)}$ then $\{z\}+J\cdot T_{(\Si,z)}\sseteq T_{(Gz,z)}$.
\\2. If $\{z\}+ J\cdot T_{(\Si,z)}\sseteq G^{(1)}z$,  then $J\sseteq ann\quotients{T_{(\Si,z)}}{T_{(G^{(1)}z,z)}}$.
\ecor
(Proof: By the assumption $J\cdot T_{(\Si,z)}$ is $G^{(1)}$-invariant and $G^{(1)}$ is unipotent for the filtration $\{(\sqrt{J})^i\cdot T_{(\Si,z)}\}_{i\ge0}$.
Thus one uses Theorem \ref{Thm.Intro.Linearization} for $M_i=(\sqrt{J})^i\cdot T_{(\Si,z)}$.)

\beR
In part (1') the condition $J\cdot T_{(\Si,z)}\sseteq T_{(G^{(1)}z,z)}$ is non-trivial and essential.
 As a trivial example, let $\Si=M\approx R$, a free module of rank one, and the group $G=R^\times$ acts by $z\to u\cdot z$.
 Then the tangent space is an ideal, $T_{(Gz,z)}=(z)\sset R$, and $ann\quotients{T_{(M,z)}}{T_{(Gz,z)}}=(z)$. But $w\in (z)$ does not imply $w+z\stackrel{G}{\sim}z$, e.g. not for $w=-z$. On the other hand, the biggest possible $G^{(1)}$ here is obtained for the filtration $\{\cm^i\}$ of $R$ and the tangent space is $T_{(G^{(1)}z,z)}=\cm\cdot z\sset R$.
\eeR

\

The bigger the ideal $ann\quotients{T_{(\Si,z)}}{T_{(Gz,z)}}$ is, the smaller  the order of determinacy is. To quantify this we use  the Loewy length of an ideal, $ll_R(I)\le\infty$, it denotes the minimal $N\le\infty$ such that $I\supseteq\cm^N$. (Here we assume  $I\supseteq\cm^\infty$.)
\bcor\label{Thm.Intro.Corol.ord.det.via.Loewy.length}
Suppose $\Si\sseteq M$ is a free direct summand, i.e., $\Si\oplus\Si^\bot=M$ for a free submodule $\Si^\bot\sset\Mat$. Then
\[
ll_R\Big(ann\quotients{T_{(\Si,z)}}{T_{(Gz,z)}}\Big)-1\le ord^\Si_G(z)\le
 ll_R\Big(ann\quotients{T_{(\Si,z)}}{T_{(G^{(1)}z,z)}}\Big)-1.
\]
\ecor
\noindent Here $G^{(1)}\sseteq G$ is the unipotent subgroup for the filtration $\{\cm^i\cdot\Si\}_i$. By corollary \ref{Thm.Group.Actions.Tangent.Space.Unipotent.Group}:
\beq
\cm\cdot T_{(Gz,z)}\sseteq T_{(G^{(1)}z,z)}\sseteq T_{(Gz,z)}.
\eeq
 Therefore the upper/lower bounds of corollary \ref{Thm.Intro.Corol.ord.det.via.Loewy.length} differ at most by 1.

\subsubsection{}
The following consequence is stated in the geometric language, using the points of the punctured neighborhood of the origin, $Spec(R)\smin\{0\}$. We consider the tangent spaces $T_{(Gz,z)}\sseteq T_{(\Si,z)}$ as  sheaves on $Spec(R)$. Accordingly, for any point $pt\in Spec(R)$ we take the fibres  $T_{(Gz,z)}|_{pt}\sseteq T_{(\Si,z)}|_{pt}$.
\bcor\label{Thm.Intro.Corol.2}
Fix an action $G\circlearrowright M$, where $G$ is \kpd\ and of Lie type.
Suppose $R$ is Noetherian and has the relevant approximation property.
Suppose $\Si\sseteq M$ is a free direct summand.
Then $z$ is finitely-$(\Si,G)$-determined \iff for any point $pt\in Spec(R)\smin\{0\}$ holds:
 $T_{(Gz,z)}|_{pt}=T_{(\Si,z)}|_{pt}$.
\ecor
\noindent(Indeed, by  corollary \ref{Thm.Intro.Corol.ord.det.via.Loewy.length}, $z$ is finitely-$(\Si,G)$-determined \iff the annihilator of $\quotients{T_{(\Si,z)}}{T_{(Gz,z)}}$ contains $\cm^N$ for some $N<\infty$. But this means that the module $\quotients{T_{(\Si,z)}}{T_{(Gz,z)}}$ vanishes off the origin.)

\beR
If  $Spec(R)\smin\{0\}$ is smooth then  $T_{(\Si,z)}|_{pt}=T_{(Gz,z)}|_{pt}$ means that $z$ is $(\Si,G)$-stable near $pt$. For $R=\k\{\ux\}$, $\k[\![\ux]\!]$ and $G$ one of the classical groups of Singularity Theory this statement is well known, e.g. \cite[Theorem 2.1]{Wall-1981}.
\eeR

\subsubsection{Matrices with $\cm$-adic filtration}
As the simplest case, consider the filtration of $M=\Mat$ by the powers of the maximal ideal, i.e., $M_i=Mat(m,n;\cm^i)$. Fix $\Si=\Mat$, then we get:
\bcor\label{Thm.Intro.Corol3} Suppose the unipotent subgroup $G^{(1)}\sset G_{lr}\rtimes Aut_\k(R)$ is \kpd\ and of Lie type. Suppose $R$ has the relevant approximation property.
\\1. $ord^\Si_{G^{(1)}}(A)=min\{k|\ Mat(m,n;\cm^{k+1})\sseteq T_{(G^{(1)}A,A)}\}$
\\2. $ord^\Si_{G^{(1)}}(A)-1\le ord^\Si_G(A)\le ord^\Si_{G^{(1)}}(A)$.
\ecor

\bex
Suppose $R$ is one of the classical rings, $\k[\![\ux]\!]$, $\k\{\ux\}$ or $C^\infty(\R^p,0)$.
Fix $\Si=M=R^{\oplus n}$ and consider the module $R^{\oplus n}$ as the space of (formal/analytic/smooth) maps from $Spec(R)$ to $(\k^n,0)$.
 For the groups $\cR$ (the right equivalence) and $\cG_{r}$ (the contact equivalence) the corollary gives the classical criterion of \cite{Mather1968},
  reproved many times, e.g. \cite{Damon1984}, \cite{Bruce.du-Plessis.Wall}.
\eex

\subsection{Admissible ideals} Quite often the ideal $ann\quotients{T_{(\Si,z)}}{T_{(Gz,z)}}$ does not contain any $\cm^k$ for $k\in\N$, thus there is no finite determinacy in the ordinary sense. Then the natural question is to find the {\em biggest ideal} $I\sset R$ such that $z$ is finitely determined for the deformations inside $\Si\cap\Big(\{z\}+I\cdot M\Big)$. Such an ideal is called ``admissible".

Suppose $R$ is Noetherian. Consider the saturation of the annihilator,  $ann^{sat}:=\suml^\infty_{i=1}\Big(ann\big(\quotients{T_{(\Si,z)}}{T_{(Gz,z)}}\big):\cm^i\Big)$,
and the radical $J=\sqrt{ann\quotients{T_{(\Si,z)}}{T_{(Gz,z)}}}$.
Then $z$ is finitely determined for deformations by $J\cdot ann^{sat}\cdot T_{(\Si,z)}$.

 This goes along the ``admissible deformations" of \cite{Siersma83} and \cite{Pellikaan88}, see section \ref{Sec.Intro.Relation.to.Singular.Theory} for other references.

\subsection{The subsequent work}\label{Sec.PartII.of.the.work}
"Theoretically" theorem \ref{Thm.Intro.Linearization} and corollary \ref{Thm.Intro.Corol.ord.det.via.Loewy.length}  ``solve" the determinacy problem: all that remains is to understand the annihilator $ann\quotients{T_{(\Si,z)}}{T_{(Gz,z)}}$.
Note that both tangent spaces are infinite-dimensional as $\k$-vector spaces, usually uncountably generated.
Considered as $R$-modules they are of high rank and in general far from being free. Thus in practice the translation of the determinacy problem to the annihilator is not yet the full/complete answer. It remains to understand the annihilator and to interpret  the condition $ann\quotients{T_{(\Si,z)}}{T_{(Gz,z)}}\supseteq\cm^N$ in terms of the particular setup.
This is similar to the transition from the theoretical
\\ $\bullet$ ``the hypersurface germ is finitely determined \iff its miniversal deformation is finite dimensional".
\\to the more practical
\\\quad $\bullet$ ``the hypersurface germ is finitely determined \iff it has at most an isolated singularity".

In  \cite{Belitski-Kerner2}, \cite{Belitski-Kerner3} we do this step, we compute (or at least bound) the annihilator for a variety of actions $G\circlearrowright \Si$.

\subsection{Relation to Singularity Theory}\label{Sec.Intro.Relation.to.Singular.Theory}
\subsubsection{}
Consider the rather particular situation:
$\k\in\R,\C$;   $R\in\k[\![\ux]\!], \k\{\ux\}, C^\infty(\R^p,0)$; $\Si=M=R^{\oplus n}$.
 Then $M$ can be considered as the space of (formal/analytic/smooth) maps from $(\k^p,0)$ to $(\k^n,0)$.
In this case instead of the group $Aut_\k(R)$ one takes the local changes of coordinates.
The two groups often coincide see section \ref{Sec.Group.Action.Aut(R).coord.changes}, and thus induces the classical right equivalence, $\cR$.
 (For a general ``non-geometric'' ring the two groups differ significantly.)

The group $GL_R(n)\rtimes Aut_\k(R)$ induces the classical contact equivalence, $\cK$. Indeed, the maps $Spec(R)\stackrel{f_1,f_2}{\to}(\k^n,0)$ are contact equivalent if $f_1(\ux)=F(f_2(\phi(\ux)),\ux)$, where $F(\ua,\ub)=L_\ub(\ua)+(h.o.t)$. Here $L_\ub(\ua)$ is linear in $\ua$, with coefficients depending on $\ub$, and $L_0$ is invertible. The $(h.o.t)$ denotes the terms at least quadratic in $\ua$.
 Thus, for a given $f_2(\ux)$ one can present: $F(f_2(\phi(\ux)),\ux)=L_{\ux,f_2,\phi}(f_2(\phi(\ux)))$, an expression linear in $f_2(\phi(\ux))$, with coefficients that depend on $\ux,f_2,\phi$.
 As $f_2(\ux)$ is fixed, write down the full dependence of $L_{\ux,f_2,\phi}$ on $\ux$ to get an element of $GL(n,R)$.
 In total, for a fixed $f_2$ we get precisely an element of $\cG_r$.

More generally, for any $R$ and $M=R^{\oplus n}$, the group actions $G\circlearrowright R$ and $G\circlearrowright(\k^n,0)$ induce the action $G\circlearrowright M$ of a subgroup $G\sseteq Aut_k(R)\rtimes GL_R(M)$.

\subsubsection{} The idea of finite determinacy begins with a simple observation: if $R$ is one of $\k[[\ux]]$, $\C\{\ux\}$, $C^\infty(\R^p,0)$ then many elements $f\in R$ are determined (up to the change of coordinates) by a few low order  monomials and do not depend on the higher order terms. The thorough investigations have (probably) began with the works of H.Whitney, R.Thom, B.Malgrange, J.N.Mather, J.C.Tougeron, V.I.Arnol'd and were continued by many others.
 (See \cite[\S III.2.2, pg.166]{AGLV}, \cite{AGLV2}, \cite{GLS} and \cite{Wall-1981}.)

\li The determinacy of maps for some subgroups of the contact group, $\cK$, was considered  in \cite{Gervais}.
In our notations he studied the actions $G\rtimes Aut(R)\circlearrowright R^{\oplus n}$, where $R=C^\infty(\R^p,0)$,
while $G\sseteq GL(n,R)$ is a Lie subgroup.
He proved that $z$ is $k$-$G$-determined \iff $T_{(Gz,z)}\supseteq \cm^{k+1}\cdot R^{\oplus n}$.
This was greatly extended in \cite{Damon1984} to the ``geometric subgroups" of $\cK$ and to the cases
  $R=\quotients{\k\{x_1,\dots,x_p\}}{I}$, $R=\quotients{\k[\![x_1,\dots,x_p]\!]}{I}$, $R=\quotients{C^\infty(\R^p,0)}{I}$.
  The importance of the unipotent part of the group (in our notation $G^{(1)}\sset G$) was clarified in \cite{Bruce.du-Plessis.Wall}.
\li The determinacy for functions on (singular) analytic germs has been studied in \cite{Bruce--Roberts}. (In our language this is the case of non-smooth germ $Spec(R)$ and the groups $GL_R(1)\rtimes Aut_\k(R)$, $Aut_\k(R)$ or their subgroups.)
\li Given two subgroups, $G,H\sset GL_\k(M)$, one can consider the ``$H$-equivariant subgroup", $G_H=\{g\in G|\ \forall\ h\in H:\ gh=hg\}$. The corresponding equivariant determinacy was studied for some subgroups of $\cK$, $\cR$
in \cite{Roberts-equivariant}, \cite{Wall-equivariant}.
\li The determinacy of square matrices (for $\k=\R$ or $\k=\C$, $R=\k\{x_1,\dots,x_p\}$ and $G=\cG_{lr}$)
was considered in \cite{Bruce-Tari04}, and further studied in \cite{Bruce-Goryun-Zakal02},
\cite{Bruce2003}, \cite{Goryun-Mond05}, \cite{Goryun-Zakal03}, \cite{Damon-Pike}. In particular, the generic finite
determinacy was established and the simple types were classified.
 Many results have been generalized in \cite{Cutkosky-Srinivasan}.
\li Sometimes one considers the coordinate changes that preserve a sublocus/subscheme in $Spec(R)$, i.e., an ideal of $R$. These were considered (for $C^\infty(\R^p,0)$-version) already in \cite{Mather1968},
 \cite{Siersma83}, \cite{Pellikaan88}, see also
 \cite{Pellikaan90}, \cite{Kushner92},  \cite{de Jong-van Straten1990}, \cite{de Jong-de Jong},   \cite{Siersma2000}, \cite{Grandjean}, \cite{Thilliez}, \cite{Sun-Wilson}, \cite{Brodersen}.
\li In the real-analytic case Arnol'd has initiated the study of functions on manifolds with boundaries, \cite{Arnold.Crit.Points.Mflds.Boundary}, see
also \cite{Goryunov.symmetric.sings} for the development and further references.
\li The study of finite determinacy in positive characteristic has been initiated in \cite{Pham}, \cite{Greuel-Pham}.

\subsubsection{}
Theorem \ref{Thm.Intro.Linearization} and its corollaries are linearization results. They reduce the initial question (highly non-linear in general) to the comparison of modules and computation of the annihilator.
This goes in the spirit of the classical Mather's criterion.

In many works the transition to the tangent-space level was done via the miniversal deformations, by proving that the infinitesimal versality implies versality. (In particular this restricted the scenarios to the cases where the miniversal deformation exists.) And usually the groups of equivalence were $\cR$, $\cK$ and some of their subgroups.
 Our results are more general in two ways.
\li We work with much broader class of rings, modules, group actions and filtrations.
\li We do not assume (and do not use) the existence of the miniversal deformation for the action $G\circlearrowright\Si$.

\

\section{The group actions and the tangent spaces}\label{Sec.Group.Actions}

\subsection{Automorphisms of local ring, $Aut_\k(R)$, vs the local coordinate changes, $\cR$}\label{Sec.Group.Action.Aut(R).coord.changes}
\subsubsection{} By definition an automorphism $\phi\in Aut_\k(R)$ satisfies:
\beq
\phi(a+b)=\phi(a)+\phi(b),\quad \phi(ab)=\phi(a)\phi(b),\quad \phi|_\k=Id.
\eeq
Thus $\phi$ is $\k$-linear and $\phi(\cm^q)=\cm^q$, i.e., the action of $\phi$
 is filtered, i.e.
$\phi$ is continuous in the Krull-topology.

\subsubsection{} For some rings the elements of the group $Aut_\k(R)$ are fixed by their action on the generators of $\cm$.
\bel
Fix some generators $\{x_i\}$ of $\cm$ (as an $R$-module).
Suppose $\cm^\infty=\{0\}$ and two automorphisms $\phi_1,\phi_2\in Aut_\k(R)$ satisfy: $\phi_1(x_i)=\phi_2(x_i)$ for any $i$. Then $\phi_1=\phi_2$.
\eel
\bpr
For any polynomial $p(\{x_i\})$ we have: $\phi_1(p)=\phi_2(p)$.
For any $q$, any element $f\in R$ can be presented in the form $p(\{x_i\})+\tilde{f}_q$, where $p$ is a polynomial, while $\tilde{f}_q\in\cm^q$. Therefore for any $f\in R$ we have: $\phi_1(f)-\phi_2(f)\in\capl_q\cm^q=\cm^\infty=\{0\}$.
\epr

Geometrically the lemma reads: $\phi$ is fully determined by its
action on the ``local coordinates", $\{x_i\}$, of $Spec(R)$. For
example, this holds for the ring $\k[[\ux]]$ and its sub-quotients.
In such cases and in this sense one can consider  $Aut_\k(R)$ as
``the group of local coordinate changes". In Singularity Theory this
group is denoted by $\cR$.

\

For ``non-geometric" rings there are many automorphisms not arising from the ``coordinate changes". \vspace{-0.2cm}
\bex
Fix two flat functions, $\tau_1,\tau_2\in C^\infty(\R^1,0)$, which are algebraically independent. Consider the ring $R=\R\{x\}[\{x^{-j}\tau_1\}_{j\in\N}, \{x^{-j}\tau_2\}_{j\in\N}]$. As $\tau_i$ are flat, the maximal ideal is generated by $x$.
Define $\phi\in Aut_\R(R)$ as identity on any converging power series and $\phi(\tau_1)=\tau_2$, $\phi(\tau_2)=\tau_1$. Extend this definition by linearity to the whole ring. We get a non-trivial automorphism that acts as identity on the ``local coordinate" $x$.
\eex

\subsection{Filtrations and completions of $M$, $End_\k(M)$ and $G\sseteq GL_\k(M)$}\label{Sec.Group.Actions.Completion.General}
Consider an $R$-module $M$ as just a $\k$-vector space and denote by $GL_\k(M)\sset End_\k(M)$ the group of all invertible $\k$-linear endomorphisms of $M$. This is the most inclusive ambient group for all the groups acting $\k$-linearly on $M$.
Note that this action is not $R$-linear, for example $GL_\k(M)$ contains the transformations induced by automorphism of the ring, $Aut_\k(R)$.
 As $M$ is usually uncountably generated over $\k$, the group $GL_\k(M)$ is huge.

The filtration $M_\bullet$ induces the filtrations of endomorphisms and automorphisms:
\beq
End^{(i)}_\k(M):=\{\phi|\ \phi(M_j)\sseteq M_{i+j},\ \forall\ j\}\sseteq End_\k(M),\quad
GL^{(i)}_\k(M)=GL_\k(M)\cap\Big\{\{\one\}+End^{(i)}_\k(M)\Big\}.
\eeq

In particular, $End^{(0)}_\k(M)$ is the space of all the endomorphisms compatible with the filtration, while
$End^{(1)}_\k(M)$ is the space of `nilpotent' endomorphism. (Note that $\phi(M_i)\sseteq M_{i+1}$ implies $\phi^k(M_i)\sseteq M_{i+k}$ but does not imply that $\phi^N=0$ for some $N$.)
Any subspace (submodule) $\La\sseteq End_\k(M)$ gets the induced filtration $\La^{(i)}:=\La\cap End^{(i)}_\k(M)$.

Similarly, $GL^{(0)}_\k(M)$ is the group of automorphisms
compatible with the filtration, while $GL^{(1)}_\k(M)$ is the unipotent subgroup of $GL_\k(M)$.
Any subgroup $G\sseteq GL_\k(M)$ gets  the induced filtration by the normal subgroups, $G^{(i)}:=G\cap GL^{(i)}_\k(M)\vartriangleleft G$. Here $G^{(1)}$ is the same as was defined in equation \eqref{Eq.def.of.G1}.
 The product/inverse operations on $G$ are continuous in this filtration topology, thus $G$ becomes a topological group.

The orbits of $G$ on $\quotients{M}{M_i}$ coincide with those of $\quotients{G}{G^{(i)}}$.

Abusing the letter $\pi_j$ we introduce the projections:
\beq
M\stackrel{\pi_j}{\to}\pi_j(M):=\quotients{M}{M_j},\quad
End_\k(M)\supseteq\La\stackrel{\pi_j}{\to} \pi_j(\La):=\quotients{\La}{\La^{(j)}},\quad
GL_\k(M)\supseteq G\stackrel{\pi_j}{\to}\pi_j(G):=\quotients{G}{G^{(i)}}.
\eeq
Here $j\le\infty$ and if $M_\infty=\{0\}$ then $\pi_\infty(M)=M$, $\pi_\infty(\La)=\La$, $\pi_\infty(G)=G$.

These filtrations/projections lead to the completion of the objects:
\beq
\hM:=\liml_{\longleftarrow}\pi_j(M),\quad \hLa:=\liml_{\longleftarrow}\pi_j(\La),\quad
\hG:=\liml_{\longleftarrow}\pi_j(G).
\eeq
We emphasize that in general $\pi_\infty(M)\neq\hM$, $\pi_\infty(\La)\neq\hLa$, $\pi_\infty(G)\neq\hG$.

\

The abuse of letter $\pi_j$ could lead to various confusions,  e.g.
\li between  $\pi_j(\La)$, the image of $\La$ under the ambient projection $End_\k(M)\stackrel{\pi_j}{\to}\pi_j(End_\k(M))$ or the image of $\La$ in $End_\k(\pi_j(M))$;
\li between the completion $\hLa$, the image of $\La$ under the ambient completion $End_\k(M)\to \widehat{End_\k(M)}$ or the image of $\La$ in $End_\k(\hM)$;
\li between the orbits of $\pi_j(G)$ and the orbits of $G$ on $\pi_j(M)$.

 The following  lemma ``justifies" such confusions.
\bel\label{Thm.Group.Actions.Completion.General}
1. For any subspace $\La\sset End^{(0)}_\k(M)$ there exist the functorial sequences of embeddings $\{\al_i\}$, $\{\be_i\}$ and isomorphisms $\{\ga_i\}$ making the following diagram commutative.
\[
\begin{picture}(0,0)(210,0)
\put(0,-2){$\pi_j\Big(End^{(0)}_\k(M)\Big)$}\put(0,-42){$End^{(0)}_\k(\pi_j M)$}\put(80,-25){$\pi_j(\La)$}

\put(25,-20){$\ga_j$}\put(68,-13){$\al_j$}\put(68,-38){$\be_j$}
\put(-10,0){\vector(-1,0){17}}\put(-10,-40){\vector(-1,0){17}}
\put(65,-22){\vector(-1,0){15}}\put(18,-22){\vector(-1,0){17}}
\put(35,-10){\begin{rotate}{-90}\scalebox{1.2}[1]{$\isom$}\end{rotate}}
\put(52,-12){\begin{rotate}{-25}\scalebox{2.4}[1]{$\hookleftarrow$}\end{rotate}}
\put(55,-38){\begin{rotate}{25}\scalebox{2.4}[1]{$\hookleftarrow$}\end{rotate}}
\end{picture}
\begin{picture}(0,0)(90,0)
\put(0,-2){$\pi_{j+1}\Big(End^{(0)}_\k(M)\Big)$}\put(0,-42){$End^{(0)}_\k(\pi_{j+1} M)$}\put(84,-25){$\pi_{j+1}(\La)$}

\put(16,-18){$\ga_{j+1}$}\put(70,-13){$\al_{j+1}$}\put(70,-38){$\be_{j+1}$}
\put(-10,0){\vector(-1,0){17}}\put(-10,-40){\vector(-1,0){17}}
\put(65,-22){\vector(-1,0){15}}\put(18,-22){\vector(-1,0){17}}
\put(35,-10){\begin{rotate}{-90}\scalebox{1.2}[1]{$\isom$}\end{rotate}}
\put(57,-12){\begin{rotate}{-25}\scalebox{2.4}[1]{$\hookleftarrow$}\end{rotate}}
\put(60,-38){\begin{rotate}{25}\scalebox{2.4}[1]{$\hookleftarrow$}\end{rotate}}
\end{picture}
\begin{picture}(0,0)(-20,0)
\put(0,0){\ldots}\put(0,-40){\ldots}\put(15,-22){\ldots}
\end{picture}
\begin{picture}(0,0)(-70,0)
\put(0,-2){$\widehat{End^{(0)}_\k(M)}$}\put(0,-42){$End^{(0)}_\k(\widehat{M})$}\put(76,-25){$\hLa$}

\put(20,-20){$\hat\ga$}\put(64,-13){$\hat\al$}\put(64,-38){$\hat\be$}
\put(-10,0){\vector(-1,0){17}}\put(-10,-40){\vector(-1,0){17}}
\put(55,-22){\vector(-1,0){15}}\put(13,-22){\vector(-1,0){17}}
\put(25,-7){\begin{rotate}{-90}\scalebox{1.2}[1]{$\isom$}\end{rotate}}
\put(48,-12){\begin{rotate}{-25}\scalebox{2.4}[1]{$\hookleftarrow$}\end{rotate}}
\put(51,-38){\begin{rotate}{25}\scalebox{2.4}[1]{$\hookleftarrow$}\end{rotate}}
\end{picture}
\begin{picture}(0,0)(-170,0)
\put(0,-2){$End^{(0)}_\k(M)$}\put(0,-42){$End^{(0)}_\k(M)$}\put(75,-25){$\La$}

\put(-10,0){\vector(-1,0){17}}\put(-10,-40){\vector(-1,0){17}}
\put(65,-22){\vector(-1,0){15}}\put(14,-22){\vector(-1,0){17}}
\put(25,-13){\begin{rotate}{-90}\scalebox{1.2}[1]{$=\!\!=$}\end{rotate}}
\put(48,-12){\begin{rotate}{-25}\scalebox{2.4}[1]{$\hookleftarrow$}\end{rotate}}
\put(51,-38){\begin{rotate}{25}\scalebox{2.4}[1]{$\hookleftarrow$}\end{rotate}}
\end{picture}
\vspace{2cm}
\]
2. For any subgroup $G\sset GL^{(0)}_\k(M)$ there exist the functorial sequences of embeddings $\{\al_i\}$, $\{\be_i\}$ and isomorphisms $\{\ga_i\}$ making the following diagram commutative.
\[
\begin{picture}(0,0)(210,0)
\put(0,-2){$\pi_j\Big(GL^{(0)}_\k(M)\Big)$}\put(0,-42){$GL^{(0)}_\k(\pi_j M)$}\put(80,-25){$\pi_j(G)$}

\put(20,-20){$\ga_j$}\put(68,-13){$\al_j$}\put(68,-38){$\be_j$}
\put(-10,0){\vector(-1,0){17}}\put(-10,-40){\vector(-1,0){17}}
\put(65,-22){\vector(-1,0){15}}\put(18,-22){\vector(-1,0){17}}
\put(30,-10){\begin{rotate}{-90}\scalebox{1.2}[1]{$\isom$}\end{rotate}}
\put(52,-12){\begin{rotate}{-25}\scalebox{2.4}[1]{$\hookleftarrow$}\end{rotate}}
\put(55,-38){\begin{rotate}{25}\scalebox{2.4}[1]{$\hookleftarrow$}\end{rotate}}
\end{picture}
\begin{picture}(0,0)(90,0)
\put(0,-2){$\pi_{j+1}\Big(GL^{(0)}_\k(M)\Big)$}\put(0,-42){$GL^{(0)}_\k(\pi_{j+1} M)$}\put(84,-25){$\pi_{j+1}(G)$}

\put(16,-18){$\ga_{j+1}$}\put(70,-13){$\al_{j+1}$}\put(70,-38){$\be_{j+1}$}
\put(-10,0){\vector(-1,0){17}}\put(-10,-40){\vector(-1,0){17}}
\put(65,-22){\vector(-1,0){15}}\put(18,-22){\vector(-1,0){17}}
\put(35,-10){\begin{rotate}{-90}\scalebox{1.2}[1]{$\isom$}\end{rotate}}
\put(57,-12){\begin{rotate}{-25}\scalebox{2.4}[1]{$\hookleftarrow$}\end{rotate}}
\put(60,-38){\begin{rotate}{25}\scalebox{2.4}[1]{$\hookleftarrow$}\end{rotate}}
\end{picture}
\begin{picture}(0,0)(-20,0)
\put(0,0){\ldots}\put(0,-40){\ldots}\put(15,-22){\ldots}
\end{picture}
\begin{picture}(0,0)(-70,0)
\put(0,-2){$\widehat{GL^{(0)}_\k(M)}$}\put(0,-42){$GL^{(0)}_\k(\widehat{M})$}\put(70,-25){$\widehat{G}$}

\put(20,-20){$\hat\ga$}\put(58,-13){$\hat\al$}\put(58,-38){$\hat\be$}
\put(-10,0){\vector(-1,0){17}}\put(-10,-40){\vector(-1,0){17}}
\put(57,-22){\vector(-1,0){15}}\put(13,-22){\vector(-1,0){17}}
\put(25,-10){\begin{rotate}{-90}\scalebox{1.2}[1]{$\isom$}\end{rotate}}
\put(42,-12){\begin{rotate}{-25}\scalebox{2.4}[1]{$\hookleftarrow$}\end{rotate}}
\put(45,-38){\begin{rotate}{25}\scalebox{2.4}[1]{$\hookleftarrow$}\end{rotate}}
\end{picture}
\begin{picture}(0,0)(-170,0)
\put(0,-2){$GL^{(0)}_\k(M)$}\put(0,-42){$GL^{(0)}_\k(M)$}\put(70,-25){$G$}

\put(-10,0){\vector(-1,0){17}}\put(-10,-40){\vector(-1,0){17}}
\put(55,-22){\vector(-1,0){15}}\put(10,-22){\vector(-1,0){17}}
\put(20,-10){\begin{rotate}{-90}\scalebox{1.2}[1]{$=\!\!=$}\end{rotate}}
\put(42,-12){\begin{rotate}{-25}\scalebox{2.4}[1]{$\hookleftarrow$}\end{rotate}}
\put(45,-38){\begin{rotate}{25}\scalebox{2.4}[1]{$\hookleftarrow$}\end{rotate}}
\end{picture}\vspace{2cm}
\]
3. For the completion map $M\stackrel{\widehat{(\ )}}{\to}\hM$ denote the image of $z\in M$ in $\hM$ by $\hz$ and the completion
  of $Gz$ by $\widehat{(Gz)}\sset\hM$.
For any $g\in GL^{(0)}_\k(M)$ and any $z\in M$ holds: $\widehat{gz}=\hg\hz$. For any $G\sseteq GL_\k(M)$ holds:
 $\widehat{(Gz)}=\widehat{\be}(\hG)\hz$.
\eel
\bpr
1. The homomorphism $\al_j$ is defined by the projection $(\phi+\La^{(j)})\to(\phi+End^{(j)}_\k(M))$.  The  homomorphism $\be_j$ is defined by: $\be_j(\phi+\La^{(j)})(z+M^{(j)})=\phi(z)+M^{(j)}$.
The  homomorphism $\ga_j$ is defined by: $\ga_j(\phi+End^{(j)}_\k(M))(z+M^{(j)})=\phi(z)+M^{(j)}$.
Note that $\al_j$, $\be_j$, $\ga_j$ do not depend on the choice of representatives, as $\La^{(j)}(M)\sseteq M_j$ and $\phi(M_j)\sseteq M_j$. Furthermore, they respect the filtration, thus their images lie in $\pi_j\big(End^{(0)}_\k(M)\big)$, $End^{(0)}_\k(\pi_j(M))$, rather than just $\pi_j End_\k(M)$, $End_\k(\pi_j(M))$.
 By the direct check:
$\ga_j\circ\al_j=\be_j$.

Now we check the claimed properties for the $j<\infty$ part of the diagram.
\bei
\item The injectivity of $\al_j$ follows from $\La^{(j)}=\La\cap End^{(j)}_\k(M)$.
\item The injectivity of $\ga_j$: \quad $\ga_j(\phi)=0$ iff $\ga_j(\phi)(M)\sseteq M^{(j)}$ iff $\phi(M)\sseteq M^{(j)}$ iff $\pi_j(\phi)=0$.
\item The surjectivity of $\ga_j$. Fix some complement, $M_j\oplus M^\bot_j=M$, as vector spaces. Fix an isomorphism $\pi_j(M)\isom M^\bot_j$, using it we identify $\pi_j(M)$ as a subspace of $M$. For any $\phi\in End_\k(\pi_j(M))$ define $\tilde\phi\in End_\k(M)$ by $\tilde\phi|_{M_j}=0$ and $\tilde\phi|_{\pi_j(M)}=\phi$. Then $\ga_j(\pi_j(\tilde\phi))=\phi$.
\item The injectivity of $\be_j$ follows from $\be_j=\ga_j\circ\al_j$.
\item The horizontal maps for $j<\infty$ are induced by $\pi_{j+1}(M)\to\pi_j(M)$. The horizontal maps at the right end of the diagram are induced by $M\to\hM$.  As the definitions of $\al_j,\be_j,\ga_j$ are  uniform in $j$ we get a commutative diagram.
\eei

\

Now we check the $j=\infty$ triangle of the diagram.
By  construction, for any $i<j$ holds: $\pi_i(\al_j)=\al_i$, $\pi_i(\be_j)=\be_i$, $\pi_i(\ga_j)=\ga_i$. Therefore the sequences $\{\al_j\}$, $\{\be_j\}$, $\{\ga_j\}$ converge and we define
$\hat{\al}:=\lim \al_j$, $\hat{\be}:=\lim \be_j$, $\hat{\ga}:=\lim \ga_j$.
More precisely, for any $\{\phi_j\}\in\hLa$, $\hat{\al}(\{\phi_j\}):=\{\al_j(\phi_j)\}$ and so on.
 (These sequences converge.)

 We check the injectivity of $\hat{\al}$. Indeed, $\hat{\al}(\{\phi_i\})=0$ means: for any $i$ exists $k_i$ such that for $j\ge k_i$ holds: $\al_j(\phi_j)\in End^{(i)}_\k(M)$. But then $\phi_j\to0$, i.e., $\{\phi_i\}=0\in\hLa$.
Similarly for $\hat{\be}$, $\hat{\ga}$.

The surjectivity of $\hat{\ga}$. Given $\hat{\phi}\in End_\k(\hM)$, we define $\tilde{\phi_j}\in \pi_j(End_\k(M))$ by
$\tilde{\phi_j}=\ga^{-1}_j\circ\phi\circ\pi_j$.  The sequence $\{\tilde{\phi_j}\}$ converges in the filtration, thus $\{\tilde{\phi_j}\}\in \widehat{End_\k(M)}$. By construction,
$\hat{\ga}(\{\tilde{\phi_j}\})=\{\ga_j\tilde{\phi_j}\}=\{\phi\circ\pi_j\}\in End_\k(\hM)$.

\

The definition of these homomorphisms does not depend on a particular choice of $\La$. Therefore, for any homomorphism $\La_1\to \La_2$ of vector subspaces of $M$ one gets the morphism of the corresponding commutative diagrams. More generally, for a morphism of any two filtered vector spaces, $M_\bullet\stackrel{f}{\to} N_\bullet$ and its  restriction to a subspace, $M_\bullet\supseteq \La\stackrel{f}{\to} f(\La)\sseteq N_\bullet$, one gets the
 morphism of the corresponding commutative diagrams. In this sense the sequences of maps $\{\al_j\}$, $\{\be_j\}$, $\{\ga_j\}$ are functorial.

2. The proof is essentially the same, replace $End^{(j)}_\k(M)$  by $GL^{(j)}_\k(M)$, $\La^{(j)}$ by $G^{(j)}$ and $+$ by $\times$.

3. To verify $\widehat{gz}=\hg\hz$ is enough to prove that the projections of both sides into $\quotients{M}{M_q}$ coincide for any $q$. Which means:
\beq
\{gz\}+ M_q\stackrel{?}{=}\{g_iz_i\}_i+ M_q,\quad\forall\ q,
\eeq
here $g_i\to \hg$ and $z_i\to \hz$.

As we take the limit, we can assume (for $i\gg0$): $g_i\in g G^{(q)}$ and $z_i\in \{z\}+ M_q$. But then the equality is obvious, as the action of $G$ is filtered.

Finally, as  $\widehat{Gz}$ is the image of $Gz$ under the completion map, we verify $\widehat{Gz}=\hG\hz$ element-wise. But for $gz\in Gz$ we have: $\widehat{gz}=\hg\hz\in\hG\hz$. Hence the statement.
\epr

\beR
If instead of $End_\k(M)$, $GL_\k(M)$ one speaks about $End_R(M)$, $GL_R(M)$ the equalities of the lemma do not hold.
The natural map $\widehat{GL_R(M)}\stackrel{\hat\ga}{\to} GL_R(\hM)$ is injective but not necessarily an isomorphism when $M$ is non-free. Take the isomorphism of the ambient groups,  $\widehat{GL_\k(M)}\stackrel{\hat\ga}{\to} GL_\k(\hM)$. By the direct check, the image of $\widehat{GL_R(M)}$ lies in $GL_R(\hM)$, thus we have the injectivity.
 For a non-surjective example, let $M=\quotients{R<s_1,s_2>}{(a_1s_1,a_2s_2)}$, where $a_1,a_2\in\cm^\infty\neq\{0\}$ and the elements $a_1,a_2$ are $R$-linearly independent. Take the filtration $\{\cm^j\cdot M\}$, then $\hM=\hR<\hat{s}_1,\hat{s}_2>\approx \hR^{\oplus 2}$. Therefore $\widehat{GL_R(M)}\approx GL_\hR(1)\times GL_\hR(1)$, but $GL_R(\hM)\approx GL_\hR(2)$.
\eeR

\subsection{\kpd\ sub-groups}\label{Sec.Group.Actions.kpd.groups}
Sometimes we forget the $R$-module structure, i.e., consider $M$ just
as a $\k$-vector space, $M_\k$. (Note that $M_\k$ is almost always
uncountably generated.) Choose a Hamel basis $\{z_\al\}$ of $M_\k$.
Then any $\k$-endomorphism $\phi\in End_\k(M)$ is presented by a
$\k$-matrix (of uncountable size),
$\phi(z_\al)=\sum_{\be}\phi_{\al\be}z_\be$. The sum here is
infinite, but only a finite number of summands are non-zero.
\bed\label{Def.kpd.groups}
 We call a subgroup $G\sseteq GL_\k(M)$ `\kpd' if it is presentable in
the form
\[G=\{\phi\in  GL_\k(M)|\ \ \uF(\{\phi_{\al\be}\})=0\},\]
 here $\uF$ is a system of polynomial equations over $\k$ (each equation is in a finite number of variables, the number of equations is usually uncountable).
\eed
\bel
Being \kpd\ does not depend on the choice of Hamel's basis of $M_\k$.
\eel
\bpr
 Any two choices of Hamel bases for $M_\k$ are related by a linear transformation, $\uz=U\uw$, $\uw=U^{-1}\uz$. Here $U$ is a $\k$-matrix of infinite (possibly uncountable) size but $U$ has in each row/column only a finite number of non-zero entries.
The change of basis implies the standard transition of the representing matrix of $\phi$, i.e., $\phi\to U\phi U^{-1}$. (Note that the product is well defined over $\k$.)
 As each equation in $\uF(\{\phi_{\al\be}\})$ contains a finite number of variables we get:
 $\uF(\{\phi_{\al\be}\})=\uF(U\{\widetilde\phi_{\al\be}\}U^{-1})=\uG(\{\widetilde\phi_{\al\be}\})$.
 The later object is again a system of polynomial equations over $\k$, each equation being in a finite number of variables.\epr

It is very difficult to work with (uncountable) Hamel's basis.
Fortunately, in many examples, the initial definition of $G\sset GL_\k(M)$ goes via some conditions of the type
\beq
\text{$\big\{F(g,a,z)=0$, for any $g\in G$, $a\in R$, $z\in M\big\}$,}
\eeq
where $F$ is some explicit expression, usually a power series or even a polynomial. Thus in these particular examples we use Hamel's bases only to verify that a group is \kpd.

A \kpd-group is defined as a subgroup of $GL_\k(M)$, in particular the action $G\circlearrowright M$ is fixed.

\

We show that the class of \kpd-groups is rich enough for our considerations.
\bel\label{Thm.Group.Actions.kpd.groups}
1. The groups  $GL_\k(M)$,  $GL_R(M)$,  $Aut_\k(R)$, $GL_R(M)\rtimes Aut_\k(R)$ are \kpd\ subgroups of $GL_\k(M)$.
\\2. Suppose $G\sseteq GL_\k(M)$ is \kpd\ and a subgroup $H\sset G$ is  defined by polynomial equations (over $\k$ or $R$). Then $H$ is \kpd.
\\3. Suppose $M$ is filtered and $G\sset GL_\k(M)$ is \kpd. Then all the subgroups $G^{(q)}\sset GL_\k(M)$ are \kpd\ for $q<\infty$.
\\4. If the groups $G\sseteq GL_R(M)$, $H\sseteq Aut_\k(R)$ are \kpd\ then the group $G\rtimes H$ is \kpd.
\\5. (Equivariant version.) Suppose a subgroup $G\sseteq GL_\k(M)$ is \kpd. Fix any subgroup $H\sseteq GL_\k(M)$ and consider $G_H=\{g\in G|\ \forall\ h\in H:\ gh=hg\}$. Then $G_H\sseteq GL_\k(M)$ is \kpd.
\eel
In particular the groups $GL^{(q)}_R(M)$,  $GL_R(M)\ltimes Aut_\k(R)$, $Aut^{(q)}_\k(R_\bullet)$ are \kpd, here
 $R_\bullet$ is a decreasing filtration of $R$ by some ideals.
\bpr
 {\bf 1.} The subgroup $GL_\k(M)\sseteq GL_\k(M)$ has no defining equations at all.
 The defining conditions of the subgroup $GL_R(M)\subset GL_\k(M)$ are: $\phi(fz)=f\phi(z)$, for any $f\in R$, $z\in M$.
  Using (any) Hamel's basis of $M$ these are written as the system $\phi(f z_\al)=f\phi(z_\al)$ for any $f\in R$ and $z_\al$. Furthermore, $f z_\al=\sum_\be f_{\al\be}z_\be$, here $\{f_{\al\be}\}$ is a $\k$-matrix of  uncountable size with finite number of non-zero entries in each row/column. We have:
\beq
f\phi(z_\al)=f\suml_\be \phi_{\al\be}z_\be=\suml_\be \phi_{\al\be}\suml_\ga f_{\be\ga}z_\ga,\quad\quad
\phi(fz_\al)=\phi(\suml_\be f_{\al\be}z_\be)=\suml_\be f_{\al\be}\suml_\ga \phi_{\be\ga}z_\ga
\eeq

Thus the defining equations of $GL_R(M)\sset GL_\k(M)$ are {\em linear}:
 $\sum_{\be} \phi_{\al\be} f_{\be\ga}=\sum_{\be} f_{\al\be} \phi_{\be\ga}$, for any $\al,\ga$.

\

The subgroup $Aut_\k(R)\sset GL_\k(R)$ is defined by the conditions $\phi(ab)=\phi(a)\phi(b)$. Again, use Hamel's basis of $R$
\beq\ber
\phi(z_\al z_\be)=\phi(\suml_\ga C_{\al\be\ga}z_\ga)=\suml_\ga C_{\al\be\ga}\suml_{\ga'}\phi_{\ga\ga'}z_{\ga'},\\
\phi(z_\al)\phi(z_\be)=\suml_{\al',\be'}\phi_{\al,\al'}\phi_{\be,\be'}z_{\al'}z_{\be'}=\suml_{\al',\be'}\phi_{\al,\al'}\phi_{\be,\be'}\suml_\ga C_{\al'\be'\ga}z_{\ga}.
\eer\eeq
Here $\{C_{\al'\be'\ga}\}$ are the `structure constants' of Hamel's basis.

Thus we get  the {\em quadratic} equations:
\beq
\forall\ \al, \be, \ga:\quad
\suml_{\ga'} C_{\al\be\ga'}\phi_{\ga'\ga}=\suml_{\al',\be'}\phi_{\al,\al'}\phi_{\be,\be'} C_{\al'\be'\ga}.
\eeq
The subgroup $GL_R(M)\rtimes Aut_\k(R)\sset GL_\k(M)$ is defined by the conditions:
$\phi(f\cdot z)=\phi(f)\phi(z)$, for any $f\in R$ and $z\in M$. Fix some Hamel bases, $\{f_\al\}$ of $R$ and $\{z_\be\}$ of $M$. Then the conditions are written in the form:
$\phi(f_\al z_\be)=\phi(f_\al)\phi(z_\be)$. Use the $R$-module structure: $f_\al\cdot z_\be=\sum_\ga c_{\al\be\ga}z_\ga$, where $c_{\al\be\ga}\in\k$ are the structure constants. Again all the defining equations are quadratic.

{\bf 2.} As $G\sseteq GL_\k(M)$ is \kpd\ we only need to check that the additional equations for $H\sset G$ are polynomial when written in Hamel's basis, $\phi(z_\al)=\sum\phi_{\al\be}z_\be$. But this is immediate as a polynomial condition $p(\phi)=0$ means: $p(\phi)(z_\al)=0$ for any $\al$. And the later translates into the polynomial equations $p(\{\phi_{\al\be}\})=0$.

\

{\bf 3.} We should prove that the additional conditions of $G^{(q)}\sset G$ translate into polynomial equations. The additional conditions are:
\beq
\text{for any }\ i\in\N, g\in G:\quad (g-\one)(M_i)\sseteq M_{i+q}.
\eeq
But for each fixed $i$, these conditions are linear in $(g-\one)$. Thus for any choice of a basis of $M$ and the corresponding
representing matrix, $g(z_\al)=\sum \phi_{\al\be}z_\be$, the equations on the entries of $(g-\one)$ are linear. Thus  $G^{(q)}$ is \kpd.

\

{\bf 4.} The elements of $G\rtimes H$ are $(g,h)$. Thus, if the defining equations of $G\sseteq GL_R(M)$ are $\uG(..)=0$ and the defining equations of $H\sseteq Aut_\k(R)$ are $\uH(..)=0$, the defining equations of
 $G\rtimes H\sseteq GL_R(M)\rtimes Aut_\k(R)$ are: $\uG(g)=0$, $\uH(h)=0$. Altogether this transforms to a system of polynomial equations over $\k$.

\

5. For any Hamel basis $\{z_\al\}$ of $M$ the equations $gh(z_\al)=hg(z_\al)$ are linear in the coefficients of $g$.
\epr

\beR
Definition \ref{Def.kpd.groups} is global, it addresses the behavior of $G$ at each point. In our paper we need only some neighborhood of the
 unit element $\one\in G$. Then one speaks of a {\em locally-\kpd\ group}. Depending on the context one uses various neighborhoods, a possible version
   is $G^{(1)}$. Sometimes one needs to include also a part of $G\smin G^{(1)}$, this is done via the projection $G\stackrel{\pi_1}{\to}\pi_1(G)$.
    The group $\pi_1(G)$ is often algebraic, thus one can take some Zariski-open neighborhood $U$ of the unit element and then ask $\pi_1^{-1}(U)$ to be \kpd.
     Sometimes $\pi_1(G)$ has a natural topology, e.g. for $\k\sseteq\C$ one can take the standard topology on $\k$ and the induced topology on $\pi_1(G)$.
      (The later is then often a classical Lie group.) Then one takes $U$ in this topology and asks $\pi_1^{-1}(U)$ to be \kpd.

One checks directly that the {\em local-\kpd} version of the lemmas of this subsection holds.
\eeR

\subsection{The tangent space to a \kpd-group}\label{Sec.Group.Actions.Tangent.Space}
\bed\label{Def.kpd.group.tangent.space}
For a \kpd\ subgroup $G=\{g|\ \uF(g)=0\}\sseteq GL_\k(M)$ the tangent space at an element $g\in G$ is defined as
$T_{(G,g)}:=\Big\{\psi\in End_\k(M)|\ \ \uF(g+\ep\psi)=0\ mod(\ep^2)\Big\}$
\eed
To write the equations explicitly fix some basis $\{z_\al\}$ of $M$, then $g(z_\al)=\sum_\be g_{\al\be}z_\be$,  $\psi(z_\al)=\sum_\be \psi_{\al\be}z_\be$.
 The group $G$ is defined by some system of equations, $\big\{\uF(\{g_{\al\be}\})=0\big\}$.
As all the equations are polynomial we can expand
$\uF(g+\ep\psi)=\uF(g)+\ep\uF'_{g}\cdot\psi+\ep^2(\cdots)$. Thus the defining equations of $T_{(G,g)}\sseteq End_\k(M)$ are $\uF'_{g}\cdot\psi=0$.
In particular $T_{(G,g)}$ is always a $\k$-linear subspace of $End_\k(M)$.

As in the previous section one shows that  the subspace $T_{(G,g)}\sseteq End_\k(M)$ does not depend on the choice of Hamel's basis of $M$.

As $T_{(G,g)}\sseteq End_\k(M)$, the action $T_{(G,g)}\circlearrowright M$ is fixed.
\bed
The tangent space to the orbit, $Gz$, at a point $z\in M$, is $T_{(Gz,z)}:=T_{(G,\one)}(z)$.
\eed
The embedding $T_{(G,g)}\sseteq End_\k(M)$ induces the filtration,
\beq
\{T^{(q)}_{(G,g)}:=T_{(G,g)}\cap End^{(q)}_\k(M)\}.
\eeq
 Accordingly we define the projections and the completion:
\beq
\pi_q(T_{(G,g)}):=\quotient{T_{(G,g)}}{T^{(q)}_{(G,g)}},\quad
\widehat{T_{(G,g)}}:=\liml_{\longleftarrow}\quotient{T_{(G,g)}}{T^{(q)}_{(G,g)}}.
\eeq

One could define the filtration/completion of tangent spaces in different ways, but they are related:
\bel\label{Thm.Group.Actions.Tangent.Spaces.General.Properties} Let $M$ be a filtered module and $G,H\sseteq GL^{(0)}_\k(M)$ some \kpd\ subgroups.
\\1. $T_{(G\cap H,g)}=T_{(G,g)}\cap T_{(H,g)}$.
\\2. In particular $T_{(G^{(q)},g)}=T^{(q)}_{(G,g)}$ and
 thus $\pi_q(T_{(G,g)})=\quotient{T_{(G,g)}}{T_{(G^{(q)},g)}}$ and
   $\widehat{T_{(G,g)}}=\liml_{\longleftarrow}\quotient{T_{(G,g)}}{T_{(G^{(q)},g)}}$.
\eel
Here the tangent spaces  are $\k$-subspaces of $End_\k(M)$, $\pi_q(End_\k(M))$, $\widehat{End_\k(M)}$, and the equalities are taken in this sense.
\bpr
{\bf 1.} Let $I_G$, $I_H$ be the defining ideals of $G,H$, then the ideal of $G\cap H$ is $I_G+I_H$. Thus:
\[\text{
 \Big($\xi\in T_{(G,g)}\cap T_{(H,g)}$\Big) iff \Big($g+\ep\xi$ satisfies the equations of $I_G$, $I_H$ modulo $\ep^2$\Big) iff
\Big($g+\ep\xi\in T_{(G\cap H,g)}$\Big).
}\]

{\bf 2.} Note that $G^{(q)}=G\cap GL^{(q)}_\k(M)$, thus by Part 1 we have:
\beq
 T_{(G^{(q)},g)}=T_{(G,g)}\cap T_{(GL^{(q)}_\k(M),g)} =T_{(G,g)}\cap End^{(q)}_\k(M)=T^{(q)}_{(G,g)}.
\eeq
\epr
\bcor\label{Thm.Group.Actions.Tangent.Space.Unipotent.Group}
Let $G\sseteq GL_\k(M)$ be a \kpd\ subgroup and suppose the filtration of $M$ satisfies: $M_i\supseteq J^i\cdot M$,
for some fixed ideal $J\ssetneq R$. Then $T_{(G^{(i)},g)}\supseteq J^i\cdot T_{(G,g)}$.
\ecor
\bpr
By lemma \ref{Thm.Group.Actions.Tangent.Spaces.General.Properties}: $T_{(G^{(i)},g)}=T^{(i)}_{(G,g)}=T_{(G,g)}\cap End_\k^{(i)}(M)$.
And by the assumption: $End_\k^{(i)}(M)\supseteq J^i\cdot End_\k(M)$. Therefore $T_{(G,g)}\cap End_\k^{(i)}(M)\supseteq J^i\cdot T_{(G,g)}$.
\epr

\bex\label{Ex.Tangent.Spaces.Ambient.Groups}
1. $T_{(GL_\k(M),g)}=End_\k(M)$, \quad \quad $T_{(GL^{(q)}_\k(M),g)}=End^{(q)}_\k(M)=\{\phi\in End_\k(M)|\ \phi(M_i)\sseteq M_{i+q}\}$.
\\2. As $GL_R(M)$ is defined by linear equations we get:
 $T_{(GL_R(M),g)}=End_R(M)$ and  $T^{(q)}_{(GL_R(M),g)}=End^{(q)}_R(M)$.
\\3. $Aut_\k(R)=\{g\in GL_\k(R)|\ \ g(ab)=g(a)g(b)\}$, thus
\beq
T_{(Aut_\k(R),g)}=\Big\{\psi\in End_\k (R)|\ \ \psi(ab)=g(a)\psi(b)+\psi(a)g(b)\Big\}.
\eeq
 In particular $T_{(Aut_\k(R),\one)}=Der_\k(R)$, the module of all the $\k$-linear derivations of $R$.
 For a regular ring the module  $Der_\k(R)$  is generated by the first order partial derivatives $\{\di_j\}$.

 Similarly for a filtration $R_\bullet$ by the ideals $J_i\supsetneq J_{i+1}\cdots$ we have:
\beq
T_{(Aut^{(q)}_\k(R_\bullet),\one)}=Der^{(q)}_\k(R_\bullet):=\{\psi\in Der_\k(R)|\ \psi J_i\sseteq J_{i+q}\}.
\eeq
4. Suppose $G=\{\uF(\{\phi\})=0\}\sset GL_R(M)\rtimes Aut_\k(R)$. Then
 $T_{(G,\one)}=\Big\{\psi\in End_R(M)+Der_\k(R)(M)|\ \ F'|_{\one}\psi=0\Big\}$.
 And for $G^{(i)}\sset G$ we have:
$T_{(G^{(i)},\one)}=T_{(G,\one)}\cap End^{(i)}_R(M_\bullet)$.

5. Let $G\sseteq GL_R(M)$ and $H\sseteq Aut_\k(R)$ then $T_{(G\rtimes H,\one)}=T_{(G,\one)}+T_{(H,\one)}$ and $[T_{(G,\one)},T_{(H,\one)}]\sseteq T_{(G,\one)}$.
\eex

\subsection{Logarithm, exponent and alternative definition of the tangent space}\label{Sec.Group.Actions.Logarithm.Exponent}
In this subsection we assume that the filtered module $M$ is $\{M_j\}$-complete, written $\hM$.
 Take any subgroup  $G\sseteq GL_{\k}(\hM)$, not necessarily \kpd.
We assume that $G$ is complete \wrt the induced filtration $\{G^{(j)}\}$.

Take the unipotent subgroup $G^{(1)}\sseteq G$ and the space of nilpotent endomorphisms $End^{(1)}_\k(\hM)$.  Define the logarithmic map:
\beq
G^{(1)}\stackrel{ln}{\to}End^{(1)}_\k(\hM),\quad\quad
g\to ln(g):=\suml^\infty_{k=1}\frac{(1-g)^k}{k}
\eeq
Note that $(1-g)\hM_j\sseteq \hM_{j+1}$, thus the sum, though infinite, is a well defined ($\k$-linear, nilpotent)  operator on $\hM$.
As this logarithm is defined by the standard Taylor series, we get: $ln(g^ig^j)=ln(g^i)+ln(g^j)$ for any $g\in G^{(1)}$. But in general $ln(gh)\neq ln(g)+ln(h)$, as $g,h\in G^{(1)}$ do not commute. In particular, the image $ln(G^{(1)})$ might be not an additive subgroup of $End^{(1)}_\k(\hM)$.

Define the exponential map $ln(G^{(1)})\stackrel{exp}{\to}GL^{(1)}_\k(\hM)$  by
 $exp(\xi):=\one+\suml^\infty_{k=1}\frac{\xi^k}{k!}$. As $\xi$ is a nilpotent endomorphism, the sum
  (though infinite) is a well defined linear operator on $\hM$ and is invertible.
\bel
$exp\Big(ln(G^{(1)})\Big)=G^{(1)}$ and the maps $ln(G^{(1)})\underset{ln}{\stackrel{exp}{\rightleftarrows}}G^{(1)}$ are mutually inverse.
\eel
\bpr
Let $\xi\in ln(G^{(1)})$, then $\xi=ln(g)$ for some $g\in G^{(1)}$. Then $\exp(\xi)=exp(ln(g))=g\in G^{(1)}$.

The maps $ln$ and $exp$ are mutual inverses as they are defined by the classical Taylor series.
\epr

In the classical situation, finite dimensional groups over a field, $ln(G^{(1)})$ is the tangent space of $G^{(1)}$ at $\one$, in fact it is the Lie algebra of the group. This holds also for \kpd\ groups:
\bprop\label{Thm.Group.Action.Tangent.Space.Comparison.Lie.Group}
Given a complete filtered module, $\hM$, and  a complete \kpd\ subgroup, $G^{(1)}\sseteq GL_\k(\hM)$, unipotent \wrt the filtration $M_\bullet$, we have:
\\1. $ln(G^{(1)})$ is a Lie algebra, i.e., it is a $k$-vector subspace of $End_\k(\hM)$, closed under commutation.
\\2. $ln(G^{(1)})\sseteq T_{(G^{(1)},\one)}$.
\eprop
Note that while $T_{(G,\one)}$ is defined externally, via the embedding $G\sseteq GL_\k(\hM)$, the space $ln(G^{(1)})$ is a purely internal object.
\bpr
1. Let $\xi\in ln(G^{(1)})$ so that $exp(\xi)\in G^{(1)}$. Then for any $n\in \Z$ holds $exp(n\xi)\in G^{(1)}$, i.e., any defining equation of $G^{(1)}$ is satisfied by $exp(t\xi)$ for $t\in\Z$. As $G^{(1)}$ is \kpd, the equations are polynomial,
and $\k$ is of zero characteristic, we get: $exp(t\xi)$ satisfies all the equations for any $t\in\k$.
 Thus $t\xi\in ln(G^{(1)})$ for any $t\in\k$, i.e., $ln(G^{(1)})$ is closed under the $\k$-multiplication.

Let $\xi_1,\xi_2\in ln(G^{(1)})$, then $exp(t\xi_1),exp(t\xi_2)\in G^{(1)}$, here we consider $t\in\k$ as a parameter.
By Baker-Campbell-Hausdorff formula:
\beq
exp(t\xi_1)\cdot exp(t\xi_2)=exp\Big(t(\xi_1+\xi_2)+\frac{t^2}{2}[\xi_1,\xi_2]+t^3(\dots)\Big).
\eeq
(Note that the proof of this formula is purely formal, it does not use any topology on $\k$ or finite dimensionality of the space.)
 For any $t\in\k$ the infinite sum $t(\xi_1+\xi_2)+\frac{t^2}{2}[\xi_1,\xi_2]+t^3(\dots)$ is a well defined nilpotent operator in $End^{(1)}_\k(M)$. Therefore $t(\xi_1+\xi_2)+\frac{t^2}{2}[\xi_1,\xi_2]+t^3(\dots)\in ln(G^{(1)})$.
 As this holds for any $t\in\k$ and $ln(G^{(1)})$ is closed under the $\k$-multiplication, we get:
\beq
\al_t=(\xi_1+\xi_2)+\frac{t}{2}[\xi_1,\xi_2]+t^2(\dots)\in ln(G^{(1)}) \text{ for any } t\neq0.
\eeq
 Therefore $exp(\al_t)\in G^{(1)}$ for any $t\neq0$. But the defining equations of $G^{(1)}$ are polynomials, thus $exp(\al_t)\in G^{(1)}$ for any $t\in\k$.
 Then $\al_t\in ln(G^{(1)})$ for any $t\in\k$, in particular $\al_0=\xi_1+\xi_2\in ln(G^{(1)})$. Therefore $ln(G^{(1)})$ is closed under addition, i.e., is a $\k$-vector space.

 Now we get: $\al_t-\xi_1-\xi_2\in ln(G^{(1)})$ for any $t$, and in the same way as above we get $[\xi_1,\xi_2]\in ln(G^{(1)})$, i.e., $ln(G^{(1)})$ is a Lie algebra.

\

2. $\sseteq$ Take any $g=exp(\xi)\in G^{(1)}$. As $ln(G^{(1)})$ is a vector space we have: $exp(t\xi)\in G^{(1)}$ for any $t\in\k$.
 In fact, each equation of $G^{(1)}$ is satisfied identically by $exp(t\xi)$, where $t$ is a variable. Thus the expanded equation vanishes in all the orders.
 In particular, $F(\one+\ep\xi)\equiv0\ mod(\ep^2)$, hence $\xi\in T_{(G^{(1)},\one)}$.
\epr

\subsection{Groups of Lie type}\label{Sec.Group.Actions.Formally.Smooth.Groups}
Definition \ref{Def.kpd.group.tangent.space} of the tangent space  goes via the defining equations of $G$,
 therefore the traditional relation of $T_{(G,g)}$ to a small neighborhood of $g$ in $G$ is not apparent.
 When $M,G$ are non-complete we do not have the maps $exp()$, $ln()$ of \S\ref{Sec.Group.Actions.Logarithm.Exponent}. And even if $M,G$ are complete
 we should clarify the relation of $ln(G)$ to $T_{(G,\one)}$.

 We define the class of groups which are close to having these maps (we take the truncated versions of $exp$, $ln$).
 \bed\label{Def.formal.smoothness}
A unipotent, locally \kpd\ subgroup $G\sset GL_\k(M)$  is called of Lie-type if the following conditions are satisfied:
\bee[i.]
\item for any $g\in G$, $q>0$ holds: $\suml^q_{j=1}\frac{(1-g)^j}{j}\in T_{(G,\one)}+End^{(q)}_\k(M)$.
\item for any $\xi\in T_{(G,\one)}$, $q>0$ holds: $\suml^q_{j=0}\frac{\xi^j}{j!}\in G\cdot GL^{(q)}_\k(M)$.
\eee
 \eed
As we show below the class of Lie-type-groups is large enough and their tangent spaces behave well.

\

First we give a general method to check that $G$ is of Lie-type. For a filtered module $M$ and the completion $\hM$, the filtered action $G\circlearrowright M$
 induces the action $G\circlearrowright\hM$, by $g(\{z_i\})=\{g(z_i)\}$. This defines a homomorphism $G\stackrel{s}{\to}s(G)\sset GL_\k(\hM)$. Note that
  in general $s$ is non-injective and $s(G)$ does not coincide with the completion $\hG$. It is enough to check the conditions of definition
   \ref{Def.formal.smoothness} for $s(G)$ and $s(T_{(G,\one)})$.

Now, for any $s\in G^{(1)}$ the operator $ln(s(g))\in End^{(1)}_\k(\hM)$ is well defined, though does not necessarily lie in $s(T_{(G,\one)})$.
 Similarly for any $\xi\in T_{(G,\one)}$  we have $exp(s(\xi))\in GL^{(1)}_\k(M)$. By the standard properties of $exp$, $ln$ we have:
\beq
exp(-s(\xi))\cdot\suml^q_{j=0}\frac{\xi^j}{j!}\in GL^{(q+1)}_\k(\hM),\quad
 ln(s(g))-\suml^q_{j=1}\frac{(1-g)^j}{j}\in End^{(q+1)}_\k(\hM).
\eeq
Therefore instead of checking the initial conditions of the definition it is enough to verify:
\beq
\forall\ \xi\in T_{(G,\one)},\ q>0:\  exp(s(\xi)) \in s(G)\cdot GL^{(q)}_\k(\hM),\quad
\forall\ g\in G,\ q>0:\  ln(s(g))\in s(T_{(G,\one)})+End^{(q)}_\k(\hM).
\eeq

\

The following statement shows that the class of Lie-type groups is rich enough.
\bel\label{Thm.Group.Actions.Smoothness}
1. The groups $GL^{(1)}_\k(M)$,  $GL^{(1)}_R(M)$ are of Lie-type.
\\2. Suppose $R=\quotients{\k[[\ux]]}{I}$ or $R=\quotients{\k\{\ux\}}{I}$ (analytic power series).
 Take any filtration $R\supset I_1\supset I_2\supset\cdots$ by ideals satisfying $\cap I_j\sseteq\cm^\infty$.
    Then  $Aut^{(1)}_\k(R)$ is of Lie type.
\\2'. Let $R\sseteq\k[[\ux]]$ or $R\sseteq C^\infty(\R^p,0)$ be any local subring that is closed under differentiation and admits substitutions,
 i.e., for any $f(x),g(x)\in\cm$
holds: $f(g(x))\in \cm$. (For example $R=C^\infty(\R^p,0)$ or $R=\k<\ux>$, algebraic power series.)
 Take any filtration $R\supset I_1\supset I_2\supset\cdots$ by ideals satisfying $\cap I_j\sseteq\cm^\infty$.
    Then  $Aut^{(1)}_\k(R)$ is of Lie type.
\\3. If $G^{(1)}$ is of Lie-type then $G^{(p)}$ is of Lie-type for any $p\ge1$.
\\4. Suppose the subgroups $G\sseteq GL_R(M)$ and $H\sseteq Aut_\k(R)$ are of Lie-type, then $G\rtimes H\sseteq GL_R(M)\rtimes Aut_\k(R)$ is of Lie-type.
\eel
\bpr
{\bf 1.}  $GL^{(1)}_\k(M)$: here for any $g\in GL^{(1)}_\k(M)$ holds: $1-g\in End^{(1)}_\k(M)=T_{(GL^{(1)}_\k(M),\one)}$. Thus $(1-g)^j\in T_{(GL^{(1)}_\k(M),\one)}$ and
 $\suml^q_{j=1}\frac{(1-g)^j}{j}\in T_{(GL^{(1)}_\k(M),\one)}$. The second condition is checked in the same way.

 $GL^{(1)}_R(M)$: it is enough to note that if $g$ is $R$-linear then $(1-g)^j$ is $R$-linear as well. Thus
  $(1-g)^j\in End^{(1)}_R(M)=T_{(GL^{(1)}_R(M),\one)}$ and  $\suml^q_{j=1}\frac{(1-g)^j}{j}\in T_{(GL^{(1)}_R(M),\one)}$.
 The second condition is checked in the same way.

\

{\bf 2.} First we check the case of $R=\quotients{\k[[\ux]]}{I}$, we prove that in this case $T_{(Aut^{(1)}_\k(R),\one)}$, $Aut^{(1)}_\k(R)$ admit  the
 ordinary exponentials/logarithmic maps. Indeed,  $T_{(Aut^{(1)}_\k(R),\one)}=Der^{(1)}_\k(R)$,  and for each $f\in R$
  the series $\suml_{j=0}^\infty\frac{\xi^j}{j!}(f)$  converges by the completeness of $R$. (And similarly for the series $\suml_{j=1}^\infty\frac{(1-g)^j}{j}(f)$.)
 Therefore the elements $exp(\xi)=\suml_{j=0}^\infty\frac{\xi^j}{j!}\in GL^{(1)}_\k(R)$, $ln(g)=\suml_{j=1}^\infty\frac{(1-g)^j}{j}\in End^{(1)}_\k(R)$ are well defined.
 Finally, the multiplicativity of the exponential series gives  $exp(\xi)(ab)=exp(\xi)(a)\cdot exp(\xi)(b)$,
while the Leibnitz rule for the logarithmic series gives $ln(g)(ab)=ln(g)(a)\cdot b+a\cdot ln(g)(b)$.

\

For an arbitrary ring $R$ the multiplicativity of  $exp$ and the Leibnitz property of $ln$ follow formally from the definition of the series,
 the only thing to check is that
 the operators $exp(\xi)$, $ln(g)$ are well defined, i.e., act on $R$. In fact, as $exp$, $ln$ are mutually inverse, it is enough to check just the case
  of $exp(\xi)$.

\

In the analytic case, $R=\quotients{\k\{\ux\}}{I}$, we note that $Der_\k(R)=\{\xi\in Der_\k(\k\{\ux\})|\ \xi(I)\sseteq I\}$ and
 $Aut_\k(R)=\{g\in Aut_\k(\k\{\ux\})|\ g(I)=I\}$. Therefore it is enough to establish the existence of $exp$ for the regular ring, $\k\{\ux\}$.
 So, we should check: for any $\xi=\sum a_i\di_i$, with $a_i\in \k\{\ux\}$, and any analytic series $f\in \k\{\ux\}$, the expression $exp(\xi)(f)$ is analytic.
 One way to do this is to prove that the series $\suml^\infty_{j=0}\frac{\xi^j(f)}{j!}$ converges in a small ball near the origin. (Note that each $\xi^j(f)$
  is analytic and the convergent sum of analytic functions is analytic.) Now to check the convergence it is enough to prove the bound:
\begin{multline}
\text{for any nilpotent derivation $\xi$ \ there exist} \ \ep>0,\ 1>C>0,\\\text{such that for any $k$ and any $|\ux|<\ep$ holds: } |\xi^k(f)(x)|<C^k k!.
\end{multline}
And this follows by multidimensional version of Cauchy formula.

\

{\bf 2'.} For an arbitrary ring $R$ of the statement the elements $exp(\xi)(f)$, $ln(g)(f)$  might not belong to $R$,
 see example \ref{Ex.Rings.with.Aut.not.of.Lie.type}, so the
 ordinary exponent/logarithm might not exist.
 Yet the approximations of definition \ref{Def.formal.smoothness} do exist.
 Indeed, the truncations $\suml^q_{j=1}\frac{(1-g)^j}{j}$, $\suml^q_{j=0}\frac{\xi^j}{j!}$ do act on $R$.
  And (being the truncations of $exp$ and $ln$) they satisfy the multiplicativity/Leibnitz
  rule modulo $I_q$.

    In more detail, associate to $\suml^q_{j=0}\frac{\xi^j}{j!}$ the operator $g_q\in End_\k(R)$ defined by $g_q(f(\ux)):=f(g_q(\ux))$, where
   $g_q(x_i)=\suml^q_{j=0}\frac{\xi^j(x_i)}{j!}$. As $R$ admits the substitution, $g_q$ is well defined.
 Moreover, by its construction $g_q$ is additive, multiplicative, unipotent and preserves $\k$. Therefore $g_q\in Aut^{(1)}_\k(R)$
    and $g^{-1}_q\cdot \suml^q_{j=0}\frac{\xi^j}{j!}\in GL^{(q)}_\k(R)$.

    Similarly,
    associate to $\suml^q_{j=1}\frac{(1-g)^j}{j}$ the operator $\xi_q\in End_\k(R)$ defined by $\xi_q(f(x)):=\sum \xi_q(x_i)\di_i f(x)$,
     where $\xi_q(x_i)=\suml^q_{j=1}\frac{(1-g)^j(x_i)}{j}$. Then $\xi_q\in Der^{(1)}_\k(R)$ and $\xi_q-\suml^q_{j=1}\frac{(1-g)^j}{j}\in End^{(q)}_\k(R)$.

Altogether we have  $\suml^q_{j=1}\frac{(1-g)^j}{j}\in Der^{(1)}_\k(R)+End^{(q)}_\k(R)$ and
  $\suml^q_{j=0}\frac{\xi^j}{j!}\in Aut^{(1)}_\k(R)\cdot GL^{(q)}_\k(R)$, thus  $Aut^{(1)}_\k(R)$ is of Lie type.

\

{\bf 3.} By lemma \ref{Thm.Group.Actions.Tangent.Spaces.General.Properties}: $T_{(G^{(p)},\one)}=T_{(G^{(1)},\one)}\cap End^{(p)}_\k(M)$ thus for
 $\xi\in T_{(G^{(p)},\one)}$ we get:
\beq
\suml^q_{j=0}\frac{\xi^j}{j!}\in (G\cap GL^{(p)}_\k(M))\cdot GL^{(q)}_\k(M)=G^{(p)}\cdot GL^{(q)}_\k(M).
\eeq
The case of $\suml^q_{j=1}\frac{(1-g)^j}{j}$ is similar.

\

{\bf 4.} Any element of $G\rtimes H$ is presentable in the form $g\cdot h$ and $T_{(G\rtimes H,\one)}=T_{(G,\one)}+T_{(H,\one)}$.
 Moreover, in our case we have: $[G^{(q)},H]\sseteq G^{(q+1)}$. Pass to the completion $\hM$, as explained after the definition.
   Then by  Baker-Campbell-Hausdorff formula we get:
 \beq
 ln(s(gh))\in ln(\underbrace{\dots}_{\in s(G)})+ln(\underbrace{\dots}_{\in s(H)})+ End^{(q)}_\k(\hM),\quad
  exp(s(\xi_G)+s(\xi_H))\in exp(\underbrace{\dots}_{\in s(T_{(G,\one)})})exp(\underbrace{\dots}_{\in s(T_{(H,\one)})})\cdot GL^{(q)}_\k(\hM).
 \eeq
As $G,H$ are of Lie-type we get $ln(s(gh))\in s(T_{(G\rtimes h,\one)})+ End^{(q)}_\k(\hM)$ and
 $exp(s(\xi_G)+s(\xi_H))\in s(G\rtimes H)\cdot GL^{(q)}_\k(\hM)$. This proves the statement.
\epr
\bex\label{Ex.Rings.with.Aut.not.of.Lie.type}
We have checked that $Aut^{(1)}_\k(R)$ is of Lie type for rather particular types of rings. Though we believe this holds for many
 other  rings, we list below some cases  where $Aut^{(1)}_\k(R)$ is not of Lie type, or at least the proof does not seem to be straightforward.
\bei
\item (the local ring of nodal cubic) Let $f=y^2-x^2-x^3$ and $R=\quotients{\k[x,y]_{(\cm)}}{(f)}$, the quotient of the localization at the origin.
 We claim: $Aut_\k(R)=\Z_2$, acting by $y\to -y$. To see this take the completion, $\hR=\quotients{\k[[x,y]]}{(y^2-x^2-x^3)}$,
 and change the variables, $a:=y-x\sqrt{1+x}$,
 $b:=y+x\sqrt{1+x}$. Then $\hR\approx \quotients{\k[[a,b]]}{(ab)}$. Therefore $Aut_\k(\hR)\approx \Z_2\ltimes GL_\hR(1)\times GL_\hR(1)$. Here $\Z_2$ acts
 by permutation $a\leftrightarrow b$, while $GL_\hR(1)\times GL_\hR(1)$ acts by scaling, $(a,b)\to (u_1a,u_2b)$, with $u_1,u_2\in\hR$ being invertible. Any element of
  $Aut_\k(R)$ ascends to $Aut_\k(\hR)$ but no (non-trivial) scaling descends to an automorphism of $R$. For example,
 in $(x,y)$ coordinates the scaling acts by $y\to \frac{u_1+u_2}{2}y+\frac{u_1-u_2}{2}x\sqrt{1+x}$, thus  if the scaling $(a,b)\to (u_1a,u_2b)$ acts on $R$ then $u_2=u_1$.
 But then $x^2+x^3\to (\frac{u_1+u_2}{2})^2(x^2+x^3)$, which over a non-henselian ring $R$ implies $\frac{u_1+u_2}{2}=\pm1$. Thus  $Aut_\k(R)=\Z_2$.
 On the other hand, $T_{(Aut_\k(R),\one)}=Der_\k(R)\neq\{0\}$, as it contains e.g. $\xi=\di_x f\di_y-\di_yf \di_x\neq0$.
 In fact $Der_\k(R)=Der_\k(\k[x,y]_{(\cm)})(-log(f))$ and therefore (as $f$ is a free divisor) is a free $R$-module of rank one.

\item (a regular local ring with mixed formal and analytic parts) Let $R=\C[[x]]\{y\}$, then $Der_\k(R)=R<\di_x,\di_y>$. Let $\xi=y^2\di_x\in Der_\k(R)$, therefore $exp(\xi)$
  induces the automorphism of $\hR=\k[[x,y]]$: $(x,y)\to(x+y^2,y)$. However this automorphism does not descend to any self-map of $R$, as it mixes $x,y$
   and sends a formal series $f(x)\in R$ to the formal series $f(x+y^2)\not\in R$.

\item
The ring $C^\infty(\R^p,0)$ does not admit the full exponential/logarithm. This is because there exist smooth functions whose subsequent derivatives grow
 arbitrarily fast in any neighborhood of the origin, therefore the series $\suml_j \frac{\xi^j(f)}{j!}$ does  not converge in general.
\eei
\eex

For Lie-type groups the two definitions of the tangent space coincide and the tangent space behaves well under completion.
\bprop\label{Thm.Complete.Groups.are.Formaly.Smooth}
Given a filtered module $M$ and a subgroup $G\sset GL_\k(M)$ which is (unipotent, \kpd\ and) of Lie type.

1. Suppose $M$ and $G$ are complete.  Then  for any $\xi\in T_{(G,\one)}$ holds:
 $exp(\xi)=\suml_{i=0}^\infty\frac{(\ep\xi)^i}{i!}\in G$. Therefore $ln(G)=T_{(G,\one)}$.

2. If both $G$ and $\hG$ are of Lie type then  $\widehat{T_{(G,\one)}}=T_{(\hG,\one)}$.
\eprop
\bpr
{\bf 1.} By the completeness of $M$  the series $\suml_{i=0}^\infty\frac{\xi^i}{i!}$ converges to an element $exp(\xi)\in GL_\k(M)$.
 As $G$ is of Lie type, $exp(\xi)\in\capl_{q>0}\Big(G\cdot GL^{(q)}_\k(M)\Big)$. Thus, as $G$ is complete, $exp(\xi)\in G$.

Invert the exponent (take the logarithm) to get $ln(G)\supseteq T_{(G,\one)}$.
The part $ln(G)\sseteq T_{(G,\one)}$ is proved in lemma \ref{Thm.Group.Action.Tangent.Space.Comparison.Lie.Group}, thus we get  $ln(G)=T_{(G,\one)}$.

\

{\bf 2.} $\supseteq$  Let $\hat\xi\in T_{(\hG,\one)}$, by part one we get: $\hat\xi=ln(\{g_i\})$ for some Cauchy sequence $\{g_i\in G\}$.
 As $G$ is of Lie-type, the element
 $\xi_q=\suml^q_{j=1}\frac{(1-g_q)^j}{j}$ belongs to $T_{(G,\one)}+End^{(q)}_\k(M)$. Thus the sequence $\{\xi_q\}$ converges
  to an element of $\widehat{T_{(G,\one)}}$ and $lim(\xi_q)=\hat\xi$. Thus $\hat\xi\in \widehat{T_{(G,\one)}}$.

$\sseteq$
Let $\hat\xi\in \widehat{T_{(G,\one)}}$, i.e., $\hat\xi=lim(\xi_q)$. Then $g_q=\suml^q_{j=0}\frac{\xi^j_q}{j!}\in G\cdot GL^{(q)}_\k(M)$, hence
 $exp(\hat\xi)=lim(g_q)\in \hG$. Thus $\hat{\xi}=ln(exp(\hat\xi))\in T_{(\hG,\one)}$.
\epr

\subsection{The main class of examples: groups acting on matrices}\label{Sec.Groups.Acting.on.Matrices}
Suppose $M_R$ is free of rank $mn$, identify $M\isom\Mat$ and take the filtration $\{Mat(m,n;\cm^q)\}_{q\in \N}$.
Many subgroups of $GL_R(M)\rtimes Aut_\k(R)$ are naturally related to the matrices. For example:
\bei
\item $G_r:=GL(n,R)$ acts on $\Mat$ by $A\to AU$;
\item $G_l:=GL(m,R)$ acts on $\Mat$ by $A\to UA$;
\item $G_{lr}:=G_l\times G_r$ and $\cG_{lr}:=G_{lr}\rtimes Aut_\k(R)$ acts by $A\to U\phi(A)V^{-1}$, $\phi\in Aut_\k(R)$.
\eei
Below we show that $\cG_{lr}$ and its natural subgroups are \kpd\ and compute their tangent spaces.

Choose a particular $R$-basis of $M$, whose generators are the matrices $\{e_{ij}\}_{\substack{i=1,\dots,m\\j=1,\dots,n}}$
with only one non-zero entry: the $(i,j)$'th entry, which is one.
Present an element $g\in GL_R(M)$ by a $mn\times mn$ matrix, $g(e_{ij})=\sum_{\tilde{i},\tilde{j}}\La_{i\tilde{i},j\tilde{j}}e_{\tilde{i},\tilde{j}}$.
\begin{enumerate}[i.]
\item  The subgroup $G_r=GL(n,R)\sset GL_R(M)$ is defined by the conditions: ``the elements of $G_r$ act on the columns of $A$", i.e.:
\beq
G_r=\Big\{g\in GL_R(M)|\quad g(e_{ij})=\suml_{\tilde{j}} a_{j\tilde{j}}e_{i\tilde{j}},\ \forall i,j\Big\}.
\eeq
(Here $a_{j\tilde{j}}\in R$ are independent of $i$.)
As in the proof of lemma \ref{Thm.Group.Actions.Tangent.Spaces.General.Properties} we fix Hamel's basis $\{z_\al\}$ of $End_\k(M)$ and expand $g=\sum g_\al z_\al$.
Thus the conditions induce the {\em linear} equations on $\{g_\al\}$. In fact $G_r\sset GL_R(M)$ is defined by $R$-linear equations.
 And the tangent space is
\beq
 T_{(G_r,\one)}=\{U\in End_R(M)|\ U(e_{ij})=\suml_{\tilde{j}} U_{j\tilde{j}}e_{i\tilde{j}}\}\isom Mat(n,n;R).
\eeq
Thus the tangent space to $G_r$-orbit of $A\in \Mat$ is $T_{(G_rA,A)}=Span_R\{Av\}_{(v\in Mat(n,n;R)}$.

By the direct check: $\pi_q(G_r)=GL(n,\pi_q(R))$ and $\widehat{G_r}=GL(n,\hR)$, in particular both groups are \kpd. Furthermore,
\beq
\pi_q( T_{(G_r,\one)})=\pi_q(Mat(n,n;R))=Mat(n,n;\pi_q(R))=T_{(\pi_q(G_r),\one)}
\eeq
 and similarly $\widehat{G_r}=Mat(n,n;\hR)=T_{(\hG_r,\one)}$.
\item Similarly, the subgroup $G_l:=GL(m,R)\sset GL_R(M)$ is defined by the conditions
``the elements of $G_r$ act on the rows of $A$". These are $R$-linear equations and the tangent space is:
 $ T_{(G_l,\one)}= Mat(m,m;R)$.
\item The definition of $G_{lr}:=G_l\times G_r=GL(m,R)\times GL(n,R)\sset GL_R(M)$ is:
\[
G_{lr}:=\Big\{g\in GL_R(M)|\quad  g(e_{ij})=\suml_{\tilde{i},\tilde{j}} a_{i\tilde{i}}b_{j\tilde{j}}e_{\tilde{i}\tilde{j}},\ \forall i,j\Big\}.
\]
(Here $\{a_{i\tilde{i}}\}$ do not depend on $j$, while $\{b_{j\tilde{j}}\}$ do not depend on $i$.)
These conditions induce quadratic equations on $\{\La_{i\tilde{i},j\tilde{j}}\}$, in particular $G_{lr}$ is \kpd.
 (The equations are precisely those of the standard Segre embedding.) The tangent space is:
   $ T_{(G_{lr},\one)}=Mat(m,m;R)\oplus Mat(n,n;R)$.
   And the tangent space to $G_{lr}$-orbit of $A\in \Mat$ is $T_{(G_{lr}A,A)}=Span_R\{uA,Av\}_{(u,v)\in Mat(m,m;R)\times Mat(n,n;R)}$.

By the direct check: $\pi_q(G_{lr})=GL(m,\pi_q(R))\times GL(n,\pi_q(R))$ and $\widehat{G_{lr}}=GL(m,\hR)\times GL(m,\hR)$,
 in particular both groups are \kpd. Furthermore, $\pi_q( T_{(G_{lr},\one)})=Mat(m,m;\pi_q(R))\oplus Mat(n,n;\pi_q(R))=T_{(\pi_q(G_{lr}),\one)}$ and
 similarly for $\widehat{G_{lr}}$.
 
\item The action of $Aut_\k(R)$ is considered in example \ref{Ex.Tangent.Spaces.Ambient.Groups}. One gets:
 $T_{(Aut_\k(R_\bullet)A,A)}=Span_R\{\cD(A)\}_{\cD\in Der(R_\bullet)}$.

\item $\cG_{lr}:=G_{lr}\rtimes Aut_\k(R)$ acts by $A\to U\phi(A)V^{-1}$. Here the tangent space to the orbit is:
\[T_{(\cG_{lr}A,A)}=Span_R\{uA,Av,\cD(A)\}_{(u,v,\cD)\in Mat(m,m;R)\times Mat(n,n;R)\times Der(R)}.
\]
Similarly  for $\cG_l:=G_{l}\rtimes Aut_\k(R)$ and $\cG_r:=G_{r}\rtimes Aut_\k(R)$.
\item $G_{congr}:\ A\to UAU^T$, is defined by quadratic equations, $\{(U,V)|\ VU=\one\}\sset G_{lr}$.
In particular $G_{congr}$ is \kpd. Here $T_{(G_{congr}A,A)}=Span_R\{uA+Au^T\}_{u\in Mat(m,m;R)}$.
Similarly for $\cG_{congr}:=G_{congr}\rtimes Aut_\k(R)$. As in the cases above,
    $\pi_q(G_{congr})=\{(U,V)|\ VU=\one\}\sset \pi_q G_{lr}$, and similarly for $\widehat{G_{congr}}$.
    Hence the isomorphism $\pi_q(T_{(G_{congr},\one)})\isom T_{(\pi_q(G_{congr}),\one)}$ and similarly for
    $\widehat{T_{(G_{congr},\one)}}$.

\item $G_{conj}:\ A\to UAU^{-1}$ is defined by linear equations $\{(U,V)|\ V=U\}\sset G_{lr}$,  hence is \kpd. Here $T_{(G_{conj}A,A)}=Span_R\{uA-Au\}_{u\in Mat(m,m;R)}$. Similarly for $\cG_{conj}:=G_{conj}\rtimes Aut_\k(R)$.
\item Let $G^{up}_l:=GL^{up}(m,R)$ denote the group of invertible upper triangular matrices over $R$. Consider the corresponding
action of $G^{up}_{lr}$: $A\to UAV$. Then $G^{up}_l$, $G^{up}_r$, $G^{up}_{lr}$ are defined by $R$-linear equations inside $G_{lr}$. Thus
$T_{(G^{up}_{lr}A,A)}=Span_R\{uA,Av\}_{(u,v)\in Mat^{up}(m,m;R)\times Mat^{up}(n,n;R)}$. Similarly for $\cG^{up}_{l}$, $\cG^{up}_{r}$, $\cG^{up}_{lr}$.
\end{enumerate}

\subsection{$R$-module structure on $T_{(Gz,z)}$}
We often restrict the class of \kpd-groups   to groups for which the $\k$-vector subspace
$T_{(Gz,z)}\sseteq T_{(M,z)}$ is an $R$-submodule. All the examples of the introduction and section \ref{Sec.Groups.Acting.on.Matrices} are of this type.

\bex
Below we list the typical groups for which $T_{(Gz,z)}$ is {\em not} an $R$-module.

1. Identify $\k$ with its embedding into $R$ and consider the subgroup $GL(n,\k)\sset GL(n,R)$. Then $T_{(GL(n,\k),\one)}\approx Mat(n,n;\k)$. This is naturally a $\k$-vector subspace of $T_{(GL(n,R),e)}\approx Mat(n,n;R)$, but not an $R$-submodule. The same behavior occurs for $GL(n,\k\oplus \cm^j)$ for any (fixed) $j\ge2$.

2. More generally, the tangent spaces are not $R$-modules for various ``nested" problems. Let $S\sset R$ be a subring with $dim(S)<dim(R)$. One often considers the corresponding subgroups, e.g. $G_S=GL(n,S)\sset GL(n,R)$. As $S$ is not an $R$ module, the tangent space $T_{(G_SA,A)}$ can never be an $R$-module.
\eex

\section{The jet-by-jet linearization procedure}\label{Sec.Linearization.procedure.jbj}

\subsection{Determinacy for the  action on a locally filtered set}\label{Sec.Actions.on.Locally.Filtred.Sets}
Let $\Si$ be an {\em arbitrary} set, not necessarily a  group.
We assume that $\Si$ is locally filtered, i.e., it is equipped with a collection of maps to some other sets, $\{\Si\stackrel{\si_j}{\to}Y_j\}_{j\in\N}$,
 such that for any $z\in \Si$ the local neighborhoods,  $\Si_{j}(z):=\si^{-1}_j\si_j(z)$, are decreasing, $\Si_j(z)\sset\Si_{j-1}(z)$. The filtration can be infinite, not necessarily stabilizing.

If $\Si$ is a subset of a filtered abelian group, $\Si\sseteq M=M_0\supset M_1\supset\cdots$, then it is natural to take the induced filtration:
 $\Si_j(z):=\Si\cap\Big(\{z\}+M_j\Big)$.

A group action $G\circlearrowright\Si$ is called filtered if $\si_j(z)=\si_j(w)$ implies $\si_j(gz)=\si_j(gw)$ for any $z,w\in \Si$ and any $g\in G$.
 We consider only filtered group actions.

 The elements $z,w\in \Si$ are called jet-by-jet-$G$ equivalent if $\{\si_j(g_j z)=\si_j(w)\}_{j\ge1}$
 for some sequence $\{g_j\in G\}_{j\in\N}$.  Denote this equivalence by $z\stackrel{G_{j.b.j.}}{\sim}w$.
 It means that $w\in\capl^\infty_{j=1}\si^{-1}_j\si_j(Gz)=\overline{Gz}$, this intersection is the closure of the orbit  $Gz$ in the filtration topology on $\Si$.

\bed
An element $z\in\Si$ is \jbj-$k$-determined if for any $w\in\Si$:  $\si_k(z)=\si_k(w)$
 implies $z\stackrel{G_{j.b.j}}{\sim}w$.
\eed
Equivalently: $\Si_k(z)\sset\overline{Gz}$.

The $j$'th stabilizer of an element $z\in\Si$ is the subgroup $St_j(z)=\{g\in G: \si_j(gz)=\si_j(z)\}$.

\

For action on an abelian group, $G\circlearrowright M$, for any $z\in M$, one defines the variation operator map, $G\stackrel{\De_z}{\to}M$, by $\De_z(g):=gz-z$.

\bel\label{Thm.Determinacy.Loc.Filtered.Sets}
1. $z\in\Si$ is $k$-$G^{(1)}_{j.b.j}$-determined \iff  \ $\forall j\ge k$ holds: $\si_{j+1}\Si_j(z)=\si_{j+1}(St_j(z)z)$.
\\2. Suppose $\Si\sseteq M$, where $M$ is an abelian group, $G\circlearrowright M$, and the filtration on $\Si$ is induced from that on $M$. Then
$z$ is \jbj-$k$-determined \iff $\forall j\ge k$: $\si_{j+1}(\Si_{j}(z)-z)=\si_{j+1}(\De_z(St_j(z)))$.
\eel
\bpr
1.
$\Rrightarrow$ Let $\Si_k(z)\sseteq \overline{G}(z)$. Assume $j\ge k$ then $\Si_j(z)\sseteq \overline{G(z)}$. Given $u\in\Si_j(z)$,
the equation $u=gz$ is jet-by-jet solvable for $g\in G$. In particular there is $g$ satisfying $\si_{j+1}u=\si_{j+1}gz$.
Since $u\in\Si_j(z)$ we get $\si_j(u)=\si_j(z)$, implying $g\in St_j(z)$. Hence, $\si_{j+1}\Si_j(z)\sseteq\si_{j+1}(St_j(z)z)$ for $j\ge k$. The inclusion
$\si_{j+1}\Si_j(z)\supseteq\si_{j+1}(St_j(z)z)$ is obvious.

$\Lleftarrow$ Given $u\in\Si_k(z)$, there is $g_1\in G$ such that $\si_{k+1}u=\si_{k+1}g_1 z$. Hence $g^{-1}_1u\in\Si_{k+1}(z)$ and there is
 $g_2\in G$ such that $\si_{k+2}g^{-1}_1u=\si_{k+1}g_2 z$. Continuing the  induction, there is $g_i\in G$ such that $\si_{k+i}u=\si_{k+1}g_1g_2\cdots g_i z$. Which means: $u$ is jet-by-jet equivalent to $z$. Thus $z$ is jet-by-jet-$k$-determined.

2. is immediate from part 1.
\epr

\subsection{Determinacy for the action on a filtered vector space}\label{Sec.Linearization.jbj.Determinacy.Action.Filtered.Vector.Space}
Suppose now that $M_\bullet$ is a filtered abelian group, $\Si\sseteq M$, $Y_i=\Si\cap M_i$. We consider the filtration $\si_i: \Si\to Y_i$, induced by $\si_i=\pi_i|_\Si$.
\bel\label{Thm.Linearization.Lemma.1}
1. $z\in M$ is $k$-$G^{(1)}_{j.b.j.}$-determined \iff $\pi_{j+1}(M_j)=\pi_{j+1}\De_z(St_j(z)\cap G^{(1)})$ for any $j\ge k$.
\\2. The composition $St_j(z)\cap G^{(1)}\stackrel{\pi_{j+1}\De_z}{\to}\pi_{j+1}(M_j)$, $j\ge1$ is a homomorphism of groups.
\\3. In particular, if $G=G^{(1)}$ and $z$ is jet-by-jet-$k$-determined
then $\pi_{j+1}(\Si_{j}(z)-z)=\pi_{j+1}(M_j)\sseteq \pi_{j+1}(M)$ is an additive subgroup for any $j\ge k$.
\bee[i.]
\item If moreover $M$ is a filtered $\k$-vector space
 then $\pi_{j+1}\Big(\Si_{j}(z)-z\Big)=\pi_{j+1}(M_j)\sseteq \pi_{j+1}(M)$ is a vector subspace.
\item If moreover $M$ is a  filtered $R$-module 
 then $\pi_{j+1}\Big(\Si_{j}(z)-z\Big)=\pi_{j+1}(M_j)\sseteq \pi_{j+1}(M)$ is an $R$-submodule.
\eee
\eel
\bpr
1. This is just part 1 of lemma \ref{Thm.Determinacy.Loc.Filtered.Sets}.

2. First note that the image of $\De_z|_{St_j(z)}$ is indeed in $M_j$, as $\pi_j(\De_z(g))=0$ for any $g\in St_j(z)$.

Let $g\in G^{(1)}$ and $h\in St_j(z)$. Then
 $\pi_{j+1}\De_z(gh)=\pi_{j+1}\big(\De_{hz}(g)+\De_z(h)\big)=\pi_{j+1}\big(\De_{z}(g)\big)+\pi_{j+1}\big(\De_z(h)\big)$, as $\pi_j(hz)=\pi_j(z)$.

3. By part 2 of lemma \ref{Thm.Determinacy.Loc.Filtered.Sets} we have $\pi_{j+1}(\Si_{j}(z)-z)=\pi_{j+1}\De_z(St_j(z))$. Therefore, by part 1,  $\pi_{j+1}\De_z(St_j(z))=\pi_{j+1}(M_j)$.
In particular $\pi_{j+1}(\Si_{j}(z)-z)$ is an additive subgroup.

 The statements i. and ii. are now immediate.
\epr

\subsection{Properties of the exponential map}\label{Sec.Linearization.jbj.properties.exponential.map}
 Suppose the $\k$-vector space $M$ is complete \wrt the filtration $M_\bullet$ and $G^{(1)}\circlearrowright M$ is complete, as  in the assumptions of section \ref{Sec.Group.Actions.Logarithm.Exponent}. Moreover, we assume that $ln(G^{(1)})$ is a $\k$-vector space, though we do not assume that $G^{(1)}$ is \kpd.
\bel\label{Thm.Linearization.Lemma.2} Let $\xi\in ln(G^{(1)})$ and $z\in M$.
\\1.  $\pi_j(exp(\xi))\in \pi_j(St_j(z))$ \iff $\pi_j(\xi z)=0\in \quotients{M}{M_j}$.
\\2. If $\pi_j(exp(\xi))\in \pi_j(St_{j}(z))$ then $\pi_{j+1}(\De_z(exp(\xi)))=\pi_{j+1}(\xi z)$.
\eel
\bpr
1. $\Rrightarrow$
As the stabilizer is a group, $\pi_j(exp(t\xi)z)=\pi_j(z)$ for all $t\in\Z$. The left hand side of this equation is a polynomial in $t$ because $\xi$ is nilpotent. As $char(\k)=0$ and this polynomial vanishes for infinitely many (distinct) values of $t$, the equality holds for all $t\in\k$. This implies $\pi_j(\xi z)=0$.

$\Lleftarrow$ If $\pi_j(\xi z)=0$ then $\pi_j(\xi^k z)=0$, thus $\pi_j(exp(\xi))\in\pi_j(St_j(z))$.

2. The function $h(t)=\pi_{j+1}(\De_z(exp(t\xi)))$ is polynomial in $t$. By part 2 of corollary \ref{Thm.Linearization.Lemma.1}  it is additive. Thus $h(t)=tc$ where $c=h(1)=\pi_{j+1} (exp(\xi)z-z)=\pi_{j+1}(\xi z)$. (In the last equation we use part 1.)
\epr
Given $z\in M$, we define the map $ln(G^{(1)})\stackrel{T_{z,j}}{\to}\pi_j(M)$ by  $T_{z,j}(\xi)=\pi_j\xi z$.
\bcor
1. For $j\ge1$:  $\pi_{j+1}\De_z(St_{j}(z))=\pi_{j+1}\Big(ker(T_{z,j})(z)\Big)$.
\\2. $z\in M$ is $k$-$G^{(1)}_{j.b.j.}$-determined \iff $\pi_{j+1}(M_j)=\pi_{j+1}\Big(ker(T_{z,j})(z)\Big)$ for any $j\ge k$.
\ecor

\subsection{The \jbj\ linearization statement}
Consider the tangent space to the orbit $\pi_j(G^{(1)}z)$ at $\pi_j(z)$, denote it
  by $\pi_j(T_{(G^{(1)}z,z)})=\pi_j(ln(G^{(1)})z)$.
\bthe\label{Thm.Linearization.General}
In the assumptions of section \ref{Sec.Linearization.jbj.properties.exponential.map}:
$z\in M$ is $k$-$G^{(1)}_{j.b.j}$-determined \iff
$M_k$ lies in the closure of $T_{(G^{(1)}z,z)}$ in Krull's topology:
 $M_k\sseteq \overline{T_{(G^{(1)}z,z)}}=\capl_{j>0}\pi_{j}^{-1}\pi_j(T_{(G^{(1)}z,z)})$.
\ethe
\bpr
$\Lleftarrow$
By the assumption we have $\pi_{j+1}(M_j)\sseteq  T_{z,j+1}\big(ln(G^{(1)})\big)$ for $j\ge k$. For any $v\in M_j$, we have
 $\pi_{j+1}v=T_{z,j+1}(\xi)=\pi_{j+1}(\xi z)$ for some $\xi\in T_{\hG^{(1)}}$. Then $\pi_j(\xi z)=\pi_j(v)=0$.
 Therefore $\pi_{j+1}(M_j)\sseteq \pi_{j+1}\De_z(St_j(z))$, and $z$ is $k$-$G_{j.b.j.}$-determined,
 by part 1 of Lemma \ref{Thm.Linearization.Lemma.1}.

$\Rrightarrow$
Let $z$ be $k$-$G_{j.b.j.}$-determined. By  part 1 of corollary \ref{Thm.Linearization.Lemma.1} for any $w\in M_k$ and any $j>k$  the equation $\pi_j(w)=\pi_j(\xi_j z)$ is  solvable for $\xi_j\in T_{\hG^{(1)}}$. Indeed, set $\xi_j=0$, for $j\le k$, and let $\xi_i\in ln(G^{(1)})$ be such that $\pi_j(w)=\pi_j(\xi_jz)$. Then
\beq
\pi_{j+1}(\xi_jz-w)\in\pi_{j+1}(M_j)=\pi_{j+1}\Big(\De_z(St_j(z))\Big)=\pi_{j+1}\Big(ker(T_{z,j})\Big)
\eeq
Hence $\pi_{j+1}(\xi_{j+1}z)=\pi_{j+1}(w)$, with $\xi_{j+1}=\xi_j-\xi$ and $\pi_{j+1}(\xi z)=\pi_{j+1}(\xi_jz-w)$.
Thus $\pi_{j+1}(M_k)\sseteq T_{z,j+1}(ln(G^{(1)}))$.
\epr

\section{The relevant  approximation results}\label{Sec.Approximation.Properties}

\subsection{The passage from the \jbj-equivalence to the equivalence of completions}\label{Sec.Approximations.Popescu.Thm}
Given a system of equations over a local ring, $\uF(\ux,\uy)=0$ (where we denote by $\ux,\uy$ finite tuples of variables), one tries to solve  iteratively, i.e., to construct a sequence $\{\uy^{(j)}(\ux)\}$ satisfying: $F(\ux,\uy^{(j)}(\ux))\equiv0\ (mod\ \cm^j)$. If $\{\uy^{(j)}(\ux)\}$ is a Cauchy sequence,  for the filtration $\{\cm^j\}$, then its limit is a formal solution, $\hF(\ux,\hat\uy(\ux))=0$. The following fundamental result ensures a formal solution without  any Cauchy property.
\bthe \cite[Theorem 2.5]{Pfister-Popescu},  \cite[Section 3A, pg 321-355]{Popescu}
 For every $F(\ux,\uy)\in\k[\![\ux,\uy]\!]^{\oplus q}$ there exists a function $\N\stackrel{\nu}{\to}\N$ satisfying:
\[\ber
\text{if $F(\ux,\tilde\uy(\ux))\equiv0\ (mod\ \cm^{\nu(c)})$, for some $\tilde\uy(\ux)\in\cm\cdot\k[\![\ux]\!]^{\oplus p}$, then there exists}
\\\text{$\uy(\ux)\in\cm\cdot\k[\![\ux]\!]^{\oplus p}$  satisfying:} \
\text{$F(\ux,\uy(\ux))\equiv0$ and $\uy(\ux)\equiv\tilde\uy(\ux)\ (mod\ \cm^c)$.}
\eer\]
\ethe
Though this statement is for the particular ring, $\k[\![\ux,\uy]\!]$, it is easily generalized:
\bcor\label{Thm.Approximat.Pfister-Popescu.Corol}
Let $(R,\cm)$ be a local ring over $\k$, suppose the $\cm$-adic completion, $\hR$, is Noetherian. Given a map $F(\ux,\uy)\in R[\![\uy]\!]^{\oplus q}$, suppose the equation $F(\ux,\uy)=0$ has a jet-by-jet solution, i.e., there exists a sequence $\{\uy_j\in \cm^{\oplus p}\}_{j=1,..,}$ satisfying $F(\ux,\uy_j)\equiv0\ (mod\ \cm^j)$. Then there exists a formal solution:
$\hat\uy(\ux)\in\hR^{\oplus p}$, $\hat{F}(\ux,\hat\uy(\ux))=0$.
\ecor
Here the completed equation, $\hF(\ux,\uy)=0$, is defined as follows. Expand $F(x,y)=\sum a_I(\ux)\uy^I$, here $I$ is a multi-index, while $a_I(\ux)\in R$. Then $\hF(\ux,\uy)=\sum \widehat{a_I(\ux)}\hat\uy^I$, where $\widehat{a_I(\ux)}$ is the image of $a_I(\ux)$ in $\hR$.
\bpr
Suppose $R$ is regular then by Cohen's theorem $\hR=\k[\![\ux]\!]$. By the assumption the equation $\hF(\ux,\hat{\uy})=0$
 has a jet-by-jet solution. Then the Popescu theorem implies a solution over $\hR$.

In the general case, by Cohen's structure theorem, $\hR=\quotients{\k[\![\ux]\!]}{I}$ for some finitely
generated ideal $I=(f_1,\dots,f_k)$. (We assume a minimal choice of generators $f_1,..,f_k$, i.e., none of $f_i$ belongs to the ideal generated by the others.)
Take any representative $\widetilde{F}$ of $\hF$ over $\k[\![\ux,\uy]\!]$.
Then a jet-by-jet solution of $F(..)=0$ means a jet-by-jet solution of $\widetilde{F}(\ux,\uy)=\sum_i z_i f_i$ over $\k[\![\ux,\uy,\uz]\!]$.
It remains to prove that (at least for $k\gg1$)  all the components of the jet-by-jet solution $\{\uy^{(k)}(\ux),\uz^{(k)}(\ux)\}_k$
 belong to $\cm=(x_1,\dots,x_p)$. For $\uy^{(k)}(\ux)$ this holds by the initial assumption. Suppose this does not hold for $\uz^{(k)}(\ux)$, i.e., for the
 infinite amount of $k$'s: $\uz^{(k)}(0)\neq0$. Then
 (by the minimality of $f_1,..,f_k$)
 there exists a sequence $\{k_l\}$ satisfying $\uz^{(k_l)}(0)=\uv$, here $\uv$ does not depend on $k$. Then we replace $\widetilde{F}(\ux,\uy)$ by $\widetilde{F}(\ux,\uy)-\sum_i v_i f_i$ and replace $\uz^{(k)}(\ux)$ by $\uz^{(k_l)}(\ux)-\uv$.
 For the new choice of $\widetilde{F}(\ux,\uy)$, choose the refined sequence, for which
 $\uz^{(k)}(0)=0$.
 Now by  Popescu theorem there exists a formal solution $\tilde{F}(\ux,\uy(\ux))=\sum_i z_i(\ux) f_i$ over $\k[\![\ux]\!]$. Its image in $\hR=\quotients{\k[\![\ux]\!]}{I}$ is the needed formal solution of $\hF(\ux,\uy)=0$.
\epr

For any action $G\circlearrowright M$ denote the closure of the orbit by $\overline{Gz}=\capl_j \Big(Gz+\cm^j\cdot M\Big)$. We compare the image of $\overline{Gz}$ under the completion, $\widehat{\overline{Gz}}$ to $\hG\hz$.
\bcor\label{Sec.Approx.Popescu.Corol.Module.is.Closed}
Let  $M$ be a (finitely generated) $R$-module with the filtration $M_j=\cm^j\cdot M$.
Suppose the completion $\hR$ is Noetherian.
Fix a unipotent subgroup $G^{(1)}\sseteq GL_R(M)\rtimes Aut_\k(R)$. (If $G^{(1)}$ involves elements of $Aut_\k(R)$ then we assume that the action $Aut_\k(R)\circlearrowright M$ is fixed.)
Suppose the $\cm$-adic completion $\hG^{(1)}$ is defined inside $GL_\hR(\hM)\rtimes Aut_\k(\hR)$ by a system of power series equations over $\hR$.  Then for any $z\in M$ holds: $\widehat{\overline{G^{(1)}z}}=\widehat{G^{(1)}}\hz$.
\ecor
In other words, given any two elements $z',z\in M$, suppose for any $j$ holds: $z'\in G^{(1)}z+\cm^j\cdot M$. Then $\hz'\in \hG^{(1)}\hz$, i.e., there exists $\hg\in\hG^{(1)}$ satisfying: $\hz'=\hg\hz$.
\bpr
Fix $z',z\in M$ and expand $\hz=\sum \ha_i\he_i$, here $\ha_i\in \hR$ while $\{\he_i\}$ is a set of generators of $\hM$, as $\hR$-module. We are trying to resolve the system of conditions:
\beq
\Bigg\{\ber \hz'=\hU\suml_i \hat{\phi}(\ha_i)\he_i\\(\hU,\hat\phi)\in\hG.\eer
\eeq
We claim that these are power series equations in $\hU$, $\hat{\phi}$. Indeed, fix some generators $\hat\ux$ of $\hat{\cm}$ over $\hR$. Then $\ha_i$ is a power series in $\hat\ux$, hence $\hat{\phi}(\ha_i(\hat\ux))=\ha_i(\hat\phi(\hat\ux))$. Thus all our conditions are power series in the unknowns. Using the unipotence, $G^{(1)}$, we present $\hU=\one+\tilde U$
 and $\hat\phi(a)=a+\tilde{\phi}(a)$, where the entries of $\tilde U$, $\tilde{\phi}$ belong to $\cm$.

By assumption this system of power series equations has a \jbj-solution whose entries belong to $\cm$. Thus by the previous corollary we get a formal solution over $\hR$.
\epr

\subsection{The passage from the equivalence of completions to the ordinary equivalence}\label{Sec.Approximations.Artin}
To pass from the completion to the ordinary equivalence we use the following
 Artin-type approximation results.

\subsubsection{The case of linear equations}\label{Sec.Approximation.for.Lin.Eqs}
 Consider a system of linear equations: $B\uy=\uv$, where $B\in\Mat$, $\uv\in Mat(m,1;R)$, while $\uy$ is the column of indeterminates.
 Consider $B$ as the map of free $R$-modules, $R^{\oplus n}\stackrel{B}{\to}R^{\oplus m}$.
 Denote by $ann.coker(B)$ the annihilator-of-cokernel of the module,
 $ann(\quotients{R^{\oplus m}}{Im(B)})$. The following property  is standard.
\bel
Fix a system of equations  $B\uy=\uv$ over a local ring $R$. Suppose
\ls either $R$ is Noetherian
\ls or the $\cm$-completion map is surjective, $R\twoheadrightarrow\hR$, and $ann.coker(B)\supseteq\cm^\infty$.
\\If the system has a formal solution (over $\hR$) then it has an ordinary solution (over $R$).
\eel
Note that in the second part $\hR$ is not assumed Noetherian. Recall that the completion map is
surjective for the ring of germs of smooth functions, $C^\infty(\R^p,0)\twoheadrightarrow\R[\![\ux]\!]$, \cite[pg. 284, exercise 12]{Rudin-book}.
\bpr
Suppose $R$ is Noetherian.
Define  the $R$-module $M_1=Im(B)$, generated by the columns of $B$. Let $M_2=M_1+\{\uv\}$. Then $M_1\sseteq M_2$ and the
 equation has a solution over $R$ \iff $M_1=M_2$. From the exact sequence of finitely generated $R$-modules, $M_1\to M_2\to \quotients{M_2}{M_1}\to0$, we pass to $M_1\otimes \hR\to M_2\otimes \hR\to \quotients{M_2}{M_1}\otimes \hR\to0$.

As the equation is solvable over $\hR$ we have: $\quotients{M_2}{M_1}\otimes\hR=\{0\}$. But the completion of Noetherian ring is faithfully flat, \cite[Theorem 55]{Matsumura}. Thus we get:
 $\quotients{M_2}{M_1}=\{0\}$. Thus $M_1=M_2$, providing a solution over $R$.

Suppose $R$ is not necessarily Noetherian but $R\twoheadrightarrow\hR$ and the system has a formal solution. Choose its representative over $R$, say $\uy_0$, and
 look for the solution (over $R$) in the form $\uy=\uy_0+\tilde\uy$. Thus we are solving the equation $B\tilde\uy=B\uy_0-\uv$.
 Note that the right hand side is ``flat", its entries belong to $\cm^\infty$.
 But then the condition $ann.coker(B)\supseteq\cm^\infty$ implies $B\uy_0-\uv\in Im(B)$, hence the solvability over $R$.
\epr

\subsubsection{Polynomial/analytic equations and Artin approximation}\label{Sec.Approximation.Artin.Classical}
Given a local ring $(R,\cm)$, consider  a finite system of equations, $F(\ux,\uy)=0$,  with $F(\ux,\uy)\in R[\uy]^{\oplus k}$.
The ring is said to have the Artin approximation property (AP) if for any solution in the $\cm$-adic completion, $\hat\uy_0(\ux)\in \hR^{\oplus p}$, $\hF(\ux,\hat\uy_0(\ux))=0$, there exists an ordinary solution
$\uy_0(\ux)\in R^{\oplus p}$, $F(\ux,\uy_0(\ux))=0$, which can be chosen arbitrarily close to $\hat{\uy}_0(\ux)$ in the $\cm$-adic topology.
A Noetherian ring $R$ over a field of zero characteristic has AP \iff $R$ is Henselian,  \cite{K.P.P.R.M.}.

Sometimes the equations are not algebraic, e.g. this happens when the action $G\circlearrowright M$ involves a  change of coordinates. If the equations are analytic then one can use the analytic Artin approximation theorem:
\bthe\cite{Artin68}
Given  a finite set of analytic equations, $F(\ux,\uy)=0$,  over $\k\{\ux,\uy\}$,
and a formal solution, $F(\ux,\hat{\uy}(\ux))\equiv0$, there exists an analytic solution, $F(\ux,\uy(\ux))\equiv0$. Moreover, $\uy(\ux)$ can be chosen arbitrarily close to $\hat{\uy}(\ux)$ in the $\cm$-adic topology.
\ethe
While this statement is formulated for a regular ring, it holds for any analytic rings, $R_X\{\uy\}$ , where $R_X=\quotients{\k\{\ux\}}{I}$, by the same argument as in the proof of corollary \ref{Thm.Approximat.Pfister-Popescu.Corol}.

\subsubsection{Equations over $C^\infty$-rings}
Tougeron's theorem \cite{Tougeron1976} says that if an {\em analytic} equation admits a formal solution, $\hat\uy(\ux)$, then it has a $C^\infty$-solution, $\uy(\ux)$, whose Taylor expansion at the origin is precisely $\hat\uy(\ux)$.

No approximation is possible when the equations are non-analytic (because of the flat functions), even for linear equations.
\bex
Let $\tau\in\C^\infty(\R^1,0)$ be a flat function. The equation $\tau^2 y+\tau=0$ has no smooth solutions, even though its completion, the identity $0\equiv0$, has plenty of solutions.
\eex
For $C^\infty$-equations the approximation holds with some additional restrictions of {\L}ojasiewicz type,
\cite[\S5]{Belitski-Kerner.Implicit}.

\subsubsection{More general equations and the Weierstrass-systems}\label{Sec.Approximation.Weierstrass}
If the equations are neither polynomial nor analytic then the approximation statement still holds for a particular class of rings called Weierstrass-systems, denoted $\k\lceil\!\lceil\ux\rceil\!\rceil$.

We do not give the explicit (lengthy) definition of Weierstrass-systems, as we do not work with them.
Rather we note that the simplest examples of W-systems are, \cite[Example 2.18]{Rond}: the formal power series, $\k[\![\ux]\!]$; the algebraic power series, $\k<\ux>$; the analytic power series, $\k\{\ux\}$; Gevrey power series. (In the last two cases $\k$ is a normed field.)
\bthe\cite{Denef-Lipshitz}, \cite[Theorem 2.19]{Rond}
Let $\k\lceil\!\lceil\ux\rceil\!\rceil$ be a W-system over $\k$. Suppose a system of equations $F\in \k\lceil\!\lceil\ux,\uy\rceil\!\rceil^{\oplus r}$ has a formal solution: $F(\ux,\hat\uy)\equiv0$, $\uy\in(\ux)\k[\![\ux]\!]^{\oplus m}$.
 Then there exists an ordinary solution: $\uy\in(\ux)\k\lceil\!\lceil\ux\rceil\!\rceil^{\oplus m}$, $F(\ux,\uy)\equiv0$. Moreover it can be chosen arbitrarily close to $\hat{\uy}$  in the $\cm$-adic topology.
\ethe

The last statement is formulated for regular rings, but it holds also for the ring $\quotients{S}{I}$, where $S$ is a W-system, while $I\sset S$ is a finitely generated ideal. The proof goes by the same argument as for
 Corollary \ref{Thm.Approximat.Pfister-Popescu.Corol}.

\subsection{Approximation properties for groups}\label{Sec.Approximation.Propert.Groups}
In section \ref{Sec.Linearization.procedure.jbj} we prove the linearization statement at the \jbj-level. Using Popescu's theorem, section \ref{Sec.Approximations.Popescu.Thm}, this is extended to the level of completion $\hG\circlearrowright \hM$.  To obtain the statement for the initial setup, $G\circlearrowright M$
we need the following approximation property:
\beq\label{Eq.Approximation.Property.for.Groups}\ber
\text{\em Fix a subgroup $G\sseteq GL_R(M)\rtimes Aut_\k(R)$, two elements $z',z\in M$ and the equation $gz=z'$ for $g\in G$. If there}\\
\text{\em  exists a formal solution, $\hat{g}\hz=\hz'$, $\hg\in\hG$, then there exists an ordinary solution, $g\in G$ such that $gz=z'$.}
\eer\eeq

We apply the approximation properties of section \ref{Sec.Approximations.Artin} and state the immediate consequences.
\bei
\item Suppose the conditions $gz=z'$, $g\in G$ can be written as a system of $R$-linear equations on $g$.
 (This is the case, e.g. for the groups of example \ref{Ex.Intro.Typical.Groups}, $GL_R(M)$, $GL^{(q)}_R(M)$, $G_l$, $G_r$, $G_{lr}$, $G_{conj}$.)
 Then the property (\ref{Eq.Approximation.Property.for.Groups}) holds for $G$ and arbitrary Noetherian local $R$.
 In the non-Noetherian case the property holds if
 $R\twoheadrightarrow\hR$, $\hR$ is Noetherian, and $ann.coker(B)\supseteq\cm^\infty$, where $B$ is the relevant matrix of coefficients of the linear equations.
 (As explained in section \ref{Sec.Approximation.for.Lin.Eqs} this happens e.g. for  $R=C^\infty(\R^p,0)$.)
 The $\cm^\infty$-conditions are checked by formulating an appropriate Lojasiewicz-inequality.
\item Suppose the conditions $gz=z'$, $g\in G$ can be written as a system of $R$-polynomial equations on $g$.
Then the property (\ref{Eq.Approximation.Property.for.Groups}) holds for $R$ Henselian and Noetherian.
\item If $\k$ is normed, there is a convergence notion in $R$ and the conditions are $R$-analytic then the
 property (\ref{Eq.Approximation.Property.for.Groups}) holds if $R$ is Henselian and Noetherian.
\item If the conditions $gz=z'$ are not polynomial/analytic
 then the approximation  (\ref{Eq.Approximation.Property.for.Groups}) holds at least when the ring $R$ is a  W-system or a quotient of W-system.
\eei

\section{Linearization and determinacy results over $R$}\label{Sec.Linearization.Final}

In section \ref{Sec.Linearization.procedure.jbj} we have established the linearization/finite determinacy at the \jbj\ level.
In section \ref{Sec.Linearization.proof.main.theorem.over.R} we combine these results with the relevant approximation properties of section \ref{Sec.Approximation.Properties} to prove the results over $R$.

\subsection{Proof of Theorem \ref{Thm.Intro.Linearization}}\label{Sec.Linearization.proof.main.theorem.over.R}

\bpr
 {\bf 1.} {\em Step 1.} Take the completion of $M$ \wrt the filtration $M_\bullet$.
 Then we have the action of the completion, $\hG^{(1)}\circlearrowright\hM$, defined in section \ref{Sec.Group.Actions.Completion.General}.
Now we compare the tangent space of the completion to the completion of the tangent space:
\beq
T_{(\hG^{(1)}\hz,\hz)}=T_{(\hG^{(1)},\one)}\hz
\stackrel{proposition \ref{Thm.Complete.Groups.are.Formaly.Smooth}}{=\!=\!=\!=\!=\!=}
\widehat{T_{(G^{(1)},\one)}}\hz
\stackrel{lemma\ \ref{Thm.Group.Actions.Completion.General}}{=\!=\!=\!=}
\widehat{T_{(G^{(1)},\one)}z}
=\widehat{T_{(G^{(1)}z,z)}}.
\eeq

Proposition \ref{Thm.Group.Action.Tangent.Space.Comparison.Lie.Group} ensures the equality: $T_{(\hG^{(1)},\one)}=ln(\hG^{(1)})$. Thus we are in the situation of Theorem \ref{Thm.Linearization.General}:

{\em the filtered action of a complete unipotent group on a complete module, $\hG^{(1)}\circlearrowright\hM_\bullet$,
 and $\hM_j\sseteq \widehat{T_{(G^{(1)}z,z)}}=T_{(\hG^{(1)}\hz,\hz)}$.}

Thus Theorem \ref{Thm.Linearization.General} implies: $\{\hz\}+\hM_j\sseteq \hG^{(1)}\hz+\hM_q$, for any $q>0$.
 Using $\quotients{M}{M_q}=\quotients{\hM}{\hM_q}$  we rewrite this statement for the initial (non-complete) module $M$:
\beq
\text{If $M_j\sseteq T_{(G^{(1)}z,z)}$ then }\{z\}+M_j\sseteq G^{(1)}z+M_q,\quad \text{for any $q>0$}.
\eeq
As the filtration $\{M_j\}$ is essentially decreasing, we get:
\beq\label{Eq.Main.Theorem.inside.proof}
\text{If $M_j\sseteq T_{(G^{(1)}z,z)}$ then }\{z\}+M_j\sseteq G^{(1)}z+\cm^q\cdot M,\quad \text{for any $q>0$}.
\eeq

{\em Step 2.}
The completion in Step 1 was taken \wrt the filtration $\{M_j\}$. Now we consider a different completion, with respect to the filtration $\{\cm^j\cdot M\}$. Denote this completion by $\widehat{(..)}^\cm$.
Equation \eqref{Eq.Main.Theorem.inside.proof} can be written using the closure: $\{z\}+M_j\sseteq \overline{G^{(1)}z}$.
But  $\{z\}+M_j\sseteq \overline{G^{(1)}z}$ implies$\{\hz^{\cm}\}+\widehat{M_j}^{\cm}\sseteq \widehat{\overline{G^{(1)}z}}^{\cm}$. Now, by corollary, \ref{Sec.Approx.Popescu.Corol.Module.is.Closed} we have: $\widehat{\overline{G^{(1)}z}}^{\cm}=\widehat{G^{(1)}}^{\cm}\hz^{\cm}$.
Which means: for any $w\in M_j$ the equations $z+w=gz$, $g\in G^{(1)}$ have a formal solution, over $\hR$.
Now invoke the approximation property of $R$ to get a solution over $R$. Thus $\{z\}+M_j\sseteq G^{(1)}z$.

\

{\bf 2.} The embedding $\{z\}+M_j\sseteq G^{(1)}z$ implies that of completions,
\beq
\{z\}+\hM_j\sseteq \widehat{G^{(1)}z}\stackrel{lemma\ \ref{Thm.Group.Actions.Completion.General}}{=\!=\!=\!=}\hG^{(1)}\hz.
\eeq
Thus by Theorem \ref{Thm.Linearization.General} we have: $\hM_j\sseteq \overline{T_{(\hG^{(1)}z,z)}}$.
 Now using proposition \ref{Thm.Complete.Groups.are.Formaly.Smooth}  we get:
\beq\label{Eq.inside.proof.main.thm}
\text{$M_j\sseteq T_{(G^{(1)},\one)}z+M_q$ for any $q>0$.}
\eeq
 In particular, as the filtration is essentially decreasing:
 $M_j\sseteq T_{(G^{(1)},\one)}z+\cm\cdot M_j$.

As $M_j$ is a finitely generated module over $R$, and $T_{(G^{(1)},\one)}z\sset M$ is a submodule, we use Nakayama lemma (note that we do not need Noetherianity of $R$) to get: $M_j\sseteq T_{(G^{(1)}z,z)}$.
\epr

\beR\label{Ex.Tangent.Space.Dont.need.submodule}
The condition ''$T_{(G^{(1)}z,z)}\sseteq T_{(M,z)}\approx M$ is a submodule" can be weakened to:
\beq \label{Eq.inside.remark}
\text{for some $N<\infty$ the intersection $T_{(G^{(1)}z,z)}\cap M_N\sset M$ is a submodule.}
\eeq
 Indeed, equation
\eqref{Eq.inside.proof.main.thm} ensures: any $w\in M_j$ is presentable in the form $\xi+w_{\ge N}$, where $\xi\in T_{(G^{(1)}z,z)}$ while $w_{\ge N}\in M_N$. Thus instead of proving $w\in T_{(G^{(1)}z,z)}$ it is enough to prove
$w_{\ge N}\in T_{(G^{(1)}z,z)}$. Or, equivalently, $M_N\sset T_{(G^{(1)}z,z)}$. By equation \eqref{Eq.inside.proof.main.thm} we have: $M_N\sseteq M_N\cap T_{(G^{(1)}z,z)}+M_{N+1}$. Now use Nakayama lemma for
$M_N\cap T_{(G^{(1)}z,z)}$.

One is tempted to weaken condition \eqref{Eq.inside.remark} further to: ``$T_{(G^{(1)}z,z)}\cap M_\infty\sset M$ is a submodule". This does not seem sufficient because of the following potentially dangerous example. Let $\k=\R$ or $\C$; $R=\k[[\ux]]$; $M=R$, $M_j=\cm^j$, and $G\sset Aut_\k(R)$ the subgroup of locally {\em analytic} coordinate changes.
 Then $T_{(G^{(1)},\one)}\ssetneq T_{(Aut^{(1)}_\k(R),\one)}$ but
 condition  \eqref{Eq.inside.proof.main.thm} holds for $j\ge2$ and any $q<\infty$. Furthermore,  $M_\infty=\{0\}$, hence $T_{(G^{(1)}z,z)}\cap M_\infty=\{0\}$ is trivially a submodule of $M$. But $T_{(G^{(1)}z,z)}$ does not contain any $M_j$.
\eeR

\subsection{Determinacy in families}\label{Sec.Determinacy.in.Families}
Let $M$ be a finitely generated $R$-module and fix a \kpd\ subgroup $G\sseteq GL_\k(M)$. We consider one-dimensional local families of elements in $M,G$.
 In detail, let $S=\k[[t]]$  or $S=\k\{t\}$, if $\k$ is a normed field. Define $SM=S\underset{\k}{\widehat\otimes}M$, accordingly one has $GL_S(SM)$ and
  the relative version $G_S\sset GL_S(SM)$ of $G$. One can check that $G_S$ is again \kpd, the unipotence of $G$ implies that of $G_S$ and if $G$ is of Lie type
   then so is $G_S$.

\bprop
Suppose $G\sseteq GL_\k(M)$ satisfies \eqref{Eq.assumptions.of.kpd.fs.etc} and is unipotent for the filtration $\{M_i\}$. Suppose $R$ has the relevant approximation property and $M_i\sseteq T_{(Gz,z)}$.
 If $z(t)\in SM_i$ then there exists $g(t)\in G_S$ such that $z(t)=g(t)z$.
\eprop
\noindent{Geometrically: if a one-parameter deformation of $z$ belongs to $\{z\}+T_{(Gz,z)}$ `pointwise' then it lies inside the orbit $Gz$, i.e., is $G$-equivalent to a trivial family.}
\bpr
Note that $G_S\sseteq GL_\k(SM)$ is \kpd\ and $T_{(G_S,\one)}=S\otimes_R T_{(G,\one)}$.
 Define the filtration of $SM$ by $(SM)_i:=SM_i$, this filtration is essentially decreasing. Then $G_S$ acts on $(SM_i)_\bullet$ and moreover $G_S$ is unipotent \wrt this filtration. Finally, the ring $S$ has the relevant approximation property because $R$ has it. Thus $SM\sseteq T_{(G_Sz,z)}$ implies by Theorem \ref{Thm.Intro.Linearization}: $\{z\}+SM\sseteq G_Sz$. Which means: for any family $z(t)\in S M$ there exists a family $g(t)\in G_S$ satisfying: $z(t)=g(t)z$.
\epr
\beR
Note that we assumed $S$ to be Henselian, $R[[t]]$ or $R\{t\}$.
For $S$ non-Henselian   we cannot ensure the existence of $g(t)\in G_S$. For example,
   let $M=Mat(1,1;R)$, with the congruence action $A\to UAU^T$, $U\in GL_R(1)$.
For $A\in Mat(1,1;R)$ consider the family $A+t fA$, where $f\in \cm^2$. Then the rectifying element, $g(t)=\sqrt{1+tf}$, belongs to $G_{S^{hen}}$, here $S^{hen}$  is the Henselization, but $g(t)\not\in G_S$.
\eeR

\end{document}